\documentclass[a4paper,11pt]{article}

\usepackage{latexsym}
\usepackage{mathrsfs}
\usepackage{amsfonts,amsmath,amssymb}
\usepackage{indentfirst}
\usepackage{verbatim}
\usepackage[pdftex]{graphicx}

\usepackage{epstopdf}
\usepackage{booktabs}
\usepackage{geometry}
\usepackage{bm}
\usepackage{enumerate}
\usepackage{float}

\usepackage{xcolor}

\newcommand{\bx}{\mbox{\boldmath{$x$}}}

\newcommand{\bzero}{\mbox{\boldmath{$0$}}}

\newcommand{\fb}{\mbox{\boldmath{$f$}}}

\newcommand{\bu}{\mbox{\boldmath{$u$}}}

\newcommand{\bv}{\mbox{\boldmath{$v$}}}

\newcommand{\bw}{\mbox{\boldmath{$w$}}}

\newcommand{\bsigma}{\mbox{\boldmath{$\sigma$}}}

\newcommand{\bnu}{\mbox{\boldmath{$\nu$}}}

\newtheorem{theorem}{Theorem}
\newtheorem{lemma}[theorem]{Lemma}

\numberwithin{equation}{section}
\newenvironment{pfs}[1][Proof]{\noindent\textbf{#1.} }
{\ \rule{0.75em}{0.75em}\smallskip}

\usepackage{xcolor}

\begin{document}


\begin{center}
\LARGE\bf Solving quasistatic contact problems\\ using nonsmooth optimization approach
\end{center}

\smallskip
\begin{center}
\large
Michal Jureczka\footnote{Jagiellonian University in Krakow, Faculty of Mathematics and Computer Science,
Lojasiewicza 6, 30-348 Krakow, Poland. Email: {\tt michal.jureczka@uj.edu.pl}},\quad
Anna Ochal\footnote{Jagiellonian University in Krakow, Faculty of Mathematics and Computer Science,
Lojasiewicza 6, 30-348 Krakow, Poland. Email: {\tt ochal@ii.uj.edu.pl}}\quad and \quad
Weimin Han\footnote{Program in Applied Mathematical and Computational Sciences (AMCS) \&
Department of Mathematics, University of Iowa, Iowa City, IA 52242, USA.
Email: {\tt weimin-han@uiowa.edu}}
\end{center}

\normalsize

\smallskip
\begin{quote}
{\bf Abstract}.
This paper is devoted to a study of time-dependent hemivariational inequality. We prove existence and uniqueness of its solution, provide fully discrete scheme and reformulate this scheme as a series of nonsmooth optimization problems. This theory is later applied to a sample quasistatic contact problem describing a viscoelastic body in frictional contact with a foundation. This contact is governed by a nonmonotone friction law with dependence on normal component of displacement and tangential component of velocity. Finally, computational simulations are performed to illustrate obtained results.

\smallskip
{\bf Keywords}. Nonmonotone friction, optimization problem, error estimate, finite element method, numerical simulations.

{\bf AMS Classification.} 35Q74, 49J40, 65K10, 65M60, 74S05, 74M15, 74M10, 74G15
\end{quote}


\section{Introduction}

Currently various engineering applications require scientific analysis of contact phenomena. In order to study such problems a framework in the form of variational and hemivariational inequalities complemented by theory of Finite Element Method emerged. Literature in the field of Contact Mechanics grows rapidly. Main ideas and mathematical tools were introduced in monographs $\cite{C, HMP, KO, MH, P}$. Other comprehensive studies include $\cite{HS, MOS, SM, W}$. Theory of nonsmooth optimization used in our approach was presented in $\cite{BKM}$ and applied to a~sample static contact problem in $\cite{BBK}$. Numerical analysis of contact problems can be found for instance in papers $\cite{BHM, H, HSM, HSB, HSD, JO1}$ and for a recent study $\cite{HS2}$.

This paper is a continuation of our previous work $\cite{JO2}$, where we reformulated time-independent hemivariational inequality as optimization problem and used it to solve a static mechanical contact problem. Here, we consider time-dependent hemivariational inequality, and use similar idea to solve it numerically. In this case we prove existence and uniqueness of a solution using fixed point argument and apply this abstract framework to solve a~quasistatic mechanical contact problem. This approach allows us to substitute Coulomb's law of dry friction by more general, nonmonotone friction law with dependence on tangential component of velocity. Additionally, the friction bound is determined by normal component of displacement, describing penetration of the foundation.

The rest of this paper is organized as follows. In Section~$\ref{sectionHI}$ we consider an abstract problem in the form of hemivariational inequality and show that under appropriate assumptions it has a~unique solution. Section~$\ref{sectionNS}$ contains a discrete scheme which approximates solution to introduced abstract problem. We also prove the theorem concerning numerical error estimate. An application of presented theory to a~mechanical contact model is introduced in Section~$\ref{sectionACP}$. Finally, in Section~$\ref{sectionNR}$, we describe computational algorithm used to solve mechanical contact problem. Then we present results of simulations and empirical error estimation for a set of sample data.

\section{Hemivariational inequality}\label{sectionHI}

 We start with basic notation used in this paper. For a normed space $X$, we denote by $\|\cdot\|_X$ its norm, by $X^*$ its dual space and by $\langle\cdot,\cdot\rangle_{X^*\times X}$ the duality pairing of $X^*$ and $X$. By $c>0$ we denote a generic constant (value of $c$ may differ in different equations).

Let us now assume that $j\colon X \to \mathbb{R}$ is locally Lipschitz continuous. The generalized directional derivative of $j$ at $x \in X$ in the direction $v\in X$ is defined by
\begin{align*}
  &j^0(x;v) := \limsup_{y \to x, \lambda \searrow 0} \frac{j(y+\lambda v)  - j(y)}{\lambda}.
\end{align*}
The generalized subdifferential of $j$ at $x$ is a subset of the dual space $X^*$ given by
\begin{align*}
  &\partial j(x) := \{\xi \in X^*\, | \, \langle \xi, v\rangle_{X^*\times X}\leq j^0(x;v) \  \mbox{ for all } v \in X \}.
\end{align*}
If $j\colon X^n \to \mathbb{R}$ is a locally Lipschitz function of $n$ variables, then we denote by $\partial_i j$ and $j_i^0$ the Clarke subdifferential and generalized directional derivative with respect to $i$-th variable of~$j$, respectively.

Let now $V$ be a reflexive Banach space and $X$ be a Banach space. Let $\gamma \in \mathcal{L}(V, X)$ be a linear and continuous operator from $V$ to $X$, and $c_{\gamma}:=\|\gamma\|_{\mathcal{L}(V, X)}$. We denote by  $\gamma^* \colon X^* \to V^*$ the adjoint operator to $\gamma$. Let $[0,T]$ be a time interval with $T>0$.   Let $A, B \colon V \to V^*$, $J \colon X \times X \to \mathbb{R}$, $f \colon [0,T] \to V^*$ and let $u_0 \in V$. The functional $J$ is assumed to be locally Lipschitz continuous with respect to its second argument. We formulate the considered hemivariational inequality as follows.

\medskip
\noindent
\textbf{Problem $\bm{P}$:} {\it\ Find $u \colon [0,T] \to V$ such that for all\ $t \in [0,T]$,}
\[\langle  Au'(t) + B u(t), w\rangle_{V^*\times V} + J_2^0(\gamma u(t),\gamma u'(t);\gamma w )
   \geq \langle  f(t),  w\rangle_{V^*\times V}\quad\forall\,w\in V, \]
and
\[ u(0) = u_0. \]

\medskip
We introduce the following assumptions.

\medskip
\noindent
$\underline{H(A)}:$ \quad Operator $A \colon V \to V^*$ satisfies
\begin{enumerate}
\item[(a)]
     $A$ is linear and bounded,
  \item[(b)]
    $A$ is symmetric, i.e., $\langle A u, v\rangle_{V^*\times V}  = \langle A v ,  u  \rangle_{V^*\times V}$ for all $u, v \in V$,
  \item[(c)]
     there exists $m_A>0$ such that $\langle A u, u\rangle_{V^*\times V} \geq m_A \|u\|_V^2$ for all $u \in V$.
\end{enumerate}

\medskip

\noindent
$\underline{H(B)}:$ \quad Operator $B \colon V \to V^*$ is Lipschitz continuous, i.e., there exists $L_B > 0$ such that\\[2mm]
 \phantom{a} \qquad $\|Bv - Bw\|_{V^*} \leq L_B\, \|v-w\|_V \mbox{ for all } v, w \in V$. \\[-2mm]

\medskip
\noindent
$\underline{H(J)}: \quad J \colon X \times X \to \mathbb{R}$ satisfies
\begin{enumerate}
  \item[(a)]
     $J(v,\cdot)$ is locally Lipschitz continuous for all $v\in V$,
  \item[(b)]
     there exist $c_0, c_1, c_2 \geq 0$ such that $\|\partial_2 J(w, v)\|_{X^*}\le c_0+c_1\|w\|_X+c_2\|v\|_X$\\
     for all $w,  v \in X$,
  \item[(c)]
    there exist $m_{J1}, m_{J2} \geq 0$ such that\\[2mm]
    $J_2^0(w_1, v_1; v_2 - v_1) + J_2^0(w _2, v_2; v_1 - v_2)
    \leq m_{J1}\| v_1- v_2\|_X^2+ m_{J2}\| w_1 -  w_2 \|_{X} \| v_1 -  v_2\|_{X}$ \\[2mm]
    for all $ w_1,  w_2,  v_1,  v_2 \in X$.
\end{enumerate}

\medskip
\noindent
$\underline{H(f)}: \quad f \in C([0,T]; V^*)$ .

\medskip
\noindent
$\underline{H(u_0)}: \quad u_0 \in V$.

\medskip
\noindent
$\underline{(H_s)}: \quad m_A - 2m_{J1}c_\gamma^2 >0$.

\bigskip
We remark that assumption $H(J)$(c) is equivalent to the following condition
\begin{align*}
 &\langle \partial_2 J(w_1,  v_1) - \partial_2 J(w_2,  v_2), v_1 -  v_2 \rangle_{X^*\times X}
 \geq - m_{J1}\|v_1-v_2\|_X^2 - m_{J2}\| w_1- w_2\|_X \| v_1- v_2\|_X
\end{align*}
for all $w_1,w_2,v_1,v_2\in X$. In the special case when $J$ does not depend on its first variable, we obtain a relaxed monotonicity condition, i.e. for all  $v_1, v_2 \in X$
\begin{align}
 &\langle \partial J(v_1) - \partial J(v_2), v_1 -  v_2 \rangle_{X^*\times X} \geq - m_{J1} \| v_1- v_2\|_X^2.  \label{relaxed}
\end{align}
So, condition $H(J)$(c) is a generalization of $(\ref{relaxed})$.

Define an operator $K \colon L^2(0,T;V)\to C([0,T];V)$ by the formula
\[ Kv(t)= \int_0^t v(s) \, ds + u_0. \]
We can reformulate Problem $P$ as the following one.

\medskip
\noindent
\textbf{Problem $\bm{P_{hvi}}$:} {\it\ Find $v \colon [0,T] \to V$ such that for all\ $t \in [0,T]$,}
\[\langle  Av(t) + B (K v)(t),w\rangle_{V^*\times V} + J_2^0(\gamma (Kv)(t), \gamma v(t);\gamma w)
 \geq \langle  f(t),  w\rangle_{V^*\times V}\quad\forall\,w\in V. \]

Let us now present the following theorem.

\begin{theorem}\label{existence}
Under assumptions $H(A)$, $H(B)$, $H(J)$, $H(f)$, $H(u_0)$ and $(H_s)$, Problem~$P_{hvi}$ has a~unique solution $v\in C([0,T];V)$.
\end{theorem}
\begin{pfs}
We use a fixed point argument. Given $\eta \in C([0,T];V)$, define
\[ y_\eta = K\eta. \]
Then $y_\eta \in C([0,T];V)$. Consider the auxiliary problem of finding a function $v_\eta \colon [0,T] \to V$ such that
\begin{equation}
\langle  Av_\eta(t),w\rangle_{V^*\times V} + J_2^0(\gamma y_\eta(t),\gamma v_\eta(t); \gamma w)
\geq \langle f(t)-B y_\eta(t),w\rangle_{V^*\times V}
\label{2.1a}
\end{equation}
for all $w\in V$. Applying Theorem 4.2 in \cite{HS2} with $\varphi \equiv 0$, we know that there exists a unique element $v_\eta(t)\in V$
which solves this inequality for each $t\in[0,T]$. Let us show that the function $v_\eta \in C([0,T];V)$.
For simplicity, with $t_1, t_2 \in [0,T]$, denote 
$v_\eta(t_i) = v_i$, $y_\eta(t_i) = y_i$,
$f(t_i) = f_i$ for $i = 1,2$. We get that
\begin{align}
\langle Av_1,w\rangle_{V^*\times V}+J_2^0(\gamma y_1,\gamma v_1;\gamma w)&\ge\langle f_1-By_1,w\rangle_{V^*\times V},
\label{first}\\
\langle Av_2,w\rangle_{V^*\times V}+J_2^0(\gamma y_2,\gamma v_2;\gamma w)&\ge\langle f_2-By_2,w\rangle_{V^*\times V}.
\label{second}
\end{align}
Taking $w = v_2-v_1$ in (\ref{first}), $w = v_1-v_2$ in (\ref{second}) and adding, we get
\begin{align*}
\langle Av_1-Av_2,v_1-v_2\rangle_{V^*\times V}
&\le J_2^0(\gamma y_1,\gamma v_1;\gamma v_2-\gamma v_1) + J_2^0(\gamma y_2,\gamma v_2;\gamma v_1-\gamma v_2)\\
&\quad {}+ \langle  f_1 - B y_1 - f_2 + B y_2, v_1 - v_2\rangle_{V^*\times V}.
\end{align*}
Using strong monotonicity of the operator $A$ guaranteed by $H(A)$(c), assumptions  $H(B)$, $H(J)$(c) and $H(f)$ we obtain
\begin{align*}
m_A\|v_1-v_2\|_V^2 & \le m_{J1}\|\gamma v_1-\gamma v_2\|_X^2+m_{J2}\|\gamma y_1-\gamma y_2\|_X\|\gamma v_1-\gamma v_2\|_X\\
&\quad{}+ (\|f_1-f_2\|_{V^*} + L_B\|y_1 - y_2\|_V) \|v_1 - v_2\|_V\\
& \le m_{J1}c_\gamma^2 \|v_1-v_2\|_V^2+\left(\|f_1-f_2\|_{V^*}+(m_{J2}c_\gamma^2+L_B)\,\|y_1-y_2\|_V\right)\|v_1-v_2\|_V,
\end{align*}
i.e.,
\begin{equation}
\left(m_A-m_{J1}c_\gamma^2\right) \|v_1-v_2\|_V\le \|f_1-f_2\|_{V^*} +\left(m_{J2}c_\gamma^2+ L_B\right)\|y_1-y_2\|_V.
\label{2.3a}
\end{equation}
By the smallness assumption $(H_s)$, $m_A-m_{J1}c_\gamma^2>0$, and from continuity of $y_\eta$ and $f$,
we deduce that $[0,T] \ni t \mapsto v_\eta (t) \in V$
is a continuous function. This allows us to define an operator $\Lambda \colon C([0,T];V)\to C([0,T];V)$ via the relation
\[ \Lambda \eta = v_\eta \quad \mbox{for\ all\ } \eta \in C([0,T];V). \]

Let us prove that the operator $\Lambda$ has a unique fixed point $\eta^* \in C([0,T];V)$. For two arbitrary functions
$\eta_1,\eta_2\in C([0,T];V)$, let $v_i$ be the solution of \eqref{2.1a} for $\eta = \eta_i$, $i=1,2$.
Similarly to \eqref{2.3a}, we have
\[ \|v_1(t) - v_2(t)\|_V\le c\,\|y_1(t) - y_2(t)\|_V. \]
Since
\[ y_1(t) - y_2(t)= \int_0^t (\eta_1(s) - \eta_2(s)) \, ds,\]
we derive from the previous inequality that
\[ \|\Lambda \eta_1(t) - \Lambda \eta_2(t)\|_V \le c\int_0^t \| \eta_1(s) - \eta_2(s)\|_V  ds.  \]
This shows that $\Lambda$ is a history dependent operator. Applying Theorem 3.20 in \cite{HS2}, we know that the operator $\Lambda$ has a unique fixed point $\eta^* \in C([0,T]; V)$. Moreover, by the definition of $\Lambda$, $\eta^*$ is a solution to Problem~$P_{hvi}$. Uniqueness of a solution to Problem~$P_{hvi}$ is a consequence of uniqueness of fixed point.
\end{pfs}

\section{Numerical scheme}\label{sectionNS}

Let us now fix $h > 0$ and let $V^h \subset V$ be a finite dimensional space with a discretization parameter $h>0$.
For a given $N\in\mathbb{N}$, we introduce the time step $k = T/N$ and the temporal nodes $t_j = jk$, $0\le j\le N$.
Let $u_0^h \in V^h$ be an approximation of element $u_0$ in space $V^h$. We write
$v^{hk} =\{v_j^{hk}\}_{j=1}^N \subset V^h$ and define
\begin{align*}
  (K^k v^{hk})_0 = u_0^h\qquad
  (K^k v^{hk})_j = k \sum_{i=1}^{j} v_i^{hk} + u_0^h, \quad 1\le j \le N.
\end{align*}
We present the following discretized version of Problem $P_{hvi}$ in the form of an operator inclusion problem.

\medskip
\noindent
\textbf{Problem $\bm{P_{incl}^h}$:} Find $v^{hk} = \{v^{hk}_j\}_{j=1}^N \subset V^h$ such that for all $j \in \{1,\dots,N\}$
\begin{align*}
  &A v^{hk}_j + B (K^kv^{hk})_{j-1} + \gamma^* \partial_2 J \big(\gamma  (K^k v^{hk})_{j-1}, \gamma v_j^{hk} \big) \ni f_j.
\end{align*}

Now we introduce some preliminary material, namely we recall a special case of the Jensen inequality and the discrete version of Gronwall inequality. On several occasions, we will also apply the elementary inequality
\begin{equation}
\Big(\sum_{i=1}^l a_i\Big)^2 \le l\sum_{i=1}^l a_i^2,\quad a_1,\cdots,a_l \in\mathbb{R}.
\label{ele}
\end{equation}

\begin{lemma} \label{jensen} {\rm (the Jensen inequality)}
Let $I \subset \mathbb{R}$ be a set of positive measure and let $f\colon I \rightarrow \mathbb{R} $ be an integrable function. Then
\begin{equation} \nonumber
\Big(\frac{1}{|I|} \int_{I} f(s)\, ds \Big)^2 \leq  \frac{1}{|I|}
\int_{I} (f(s))^2\, ds.
\end{equation}
In particular, for $t_m, t_{m+1} \in I$, $t_m < t_{m+1}$, $t_{m+1} - t_{m} = k$, we have
\begin{equation} \label{JenIneq}
\Big(\int_{t_m}^{t_{m+1}} f(s)\, ds \Big)^2 \leq  k \int_{t_m}^{t_{m+1}} (f(s))^2\, ds.
\end{equation}
\end{lemma}

\begin{lemma} \label{gronwall} {\rm (the Gronwall inequality, \cite[Lemma 7.25]{HS})} \
Let $T$ be given. For $N>0$ we define $k=T/N$. Let $\{g_n\}_{n=1}^N$, $\{e_n\}_{n=1}^N$ be two nonnegative sequences satisfying for $c>0$ and for all $n\in\{1,\dots,N\}$
\begin{equation} \nonumber
 e_n\leq cg_n + ck\sum_{j=1}^{n-1} e_j.
\end{equation}
Then there exists a constant $\hat{c}>0$ such that
\begin{equation} \nonumber
\max_{1\leq n\leq N} e_n \leq \hat{c} \max_{1\leq n\leq N} g_n.
\end{equation}
\end{lemma}

Let us now prove the following lemma.

\begin{lemma} \label{I}
Under assumptions $H(A)$, $H(B)$, $H(J)$, $H(f)$, $H(u_0)$ and $(H_s)$ if Problem~$P_{incl}^h$ has a~solution $v^{hk}$, then it is unique and satisfies
\begin{align}
  &\| v_j^{hk} \|_V \leq c(1+\| f_j \|_{V^*}) \label{e22}
\end{align}
for all $j \in \{1,\dots,N\}$ with a positive constant $c$.
\end{lemma}
\begin{pfs}
To simplify the notation in this proof, we write $v$ instead of $v^{hk}$.
Let $v$ be a solution to Problem  $P_{incl}^h$ and let us fix $j \in \{1,\dots,N\}$. This means that there exists $\zeta_j \in  \partial_2 J \big(\gamma  (K^k v)_{j-1}, \gamma v_j \big)$ such that
\begin{align*}
 &A v_j + B (K^kv)_{j-1} + \gamma^* \zeta_j = f_j.
\end{align*}
From the definition of generalized subdifferential of $J(\gamma (K^k v)_{j-1}, \cdot)$ we have for all $w^h \in V^h$,
\begin{align*}
\langle f_j - A v_j - B (K^kv)_{j-1},  w^h\rangle_{V^*\times V}
  &= \langle \gamma^* \zeta_j,  w^h\rangle_{V^*\times V} \\
  & = \langle \zeta_j , \gamma w^h\rangle_{X^*\times X}\\
  &\leq J_2^0(\gamma (K^k v)_{j-1}, \gamma v_j;  \gamma w^h).
\end{align*}
After reformulation, we obtain discretized version of hemivariational inequality $P_{hvi}$,
\begin{equation}
\langle A v_j+B(K^kv)_{j-1},w^h\rangle_{V^*\times V}+J_2^0(\gamma (K^k v)_{j-1},\gamma v_j;\gamma w^h)
\geq \langle f_j, w^h \rangle_{V^*\times V}.
\label{9}
\end{equation}
Let us now assume that Problem $P_{incl}^h$ has two solutions $v^1$ and $v^2$. We will prove inductively
that these solutions are equal.
For $j = 1$ we get $(K^kv^1)_0 = (K^kv^2)_0 = u_0^h$. In inequality \eqref{9} for a~solution $v^1$ we set
$ w^h =v^2_1-v^1_1$, and for a solution $v^2$ we set $w^h  =  v_1^1-v^2_1$.
Then adding these inequalities, we obtain
\begin{align*}
  &\langle A v_1^1 - A v_1^2,  v_1^1 - v_1^2 \rangle_{V^*\times V} \leq J_2^0(\gamma  u_0^h, \gamma v_1^1; \gamma v_1^2 - \gamma v_1^1 ) +J_2^0(\gamma  u_0^h, \gamma v_1^2; \gamma v_1^1 - \gamma v_1^2 ).
\end{align*}
Hence, $H(A)$(c) and $H(J)$(c) yield
\begin{align*}
  &\big(m_A - m_{J1}c_\gamma^2\big) \|  v_1^1 - v_1^2\|_V^2 \leq 0,
\end{align*}
and consequently from $(H_s)$ we have $v^1_1 = v^2_1$. We now show that if $v_i^1 = v_i^2$ for $i=1,\dots,j-1$, then $v_j^1 = v_j^2$. Similarly, in~(\ref{9}) for a solution $v^1$ we set $ w^h  =  v^2_j-v^1_j$, and for a solution $v^2$ we set  $w^h  =  v_j^1-v^2_j$. Adding these inequalities we obtain
\begin{align*}
\langle A v_j^1 - A v_j^2,  v_j^1 - v_j^2 \rangle_{V^*\times V}
&\leq \langle B (K^kv^1)_{j-1} - B (K^kv^2)_{j-1},  v_j^2 - v_j^1 \rangle_{V^*\times V}\\
&\quad{}+J_2^0(\gamma (K^k v^1)_{j-1}, \gamma v_j^1; \gamma v_j^2 - \gamma v_j^1 )
+J_2^0(\gamma (K^k v^2)_{j-1}, \gamma v_j^2; \gamma v_j^1 - \gamma v_j^2 ).
\end{align*}
We observe that $(K^kv^1)_{j-1} = (K^kv^2)_{j-1}$, hence
\begin{align*}
  &\langle A v_j^1 - A v_j^2,  v_j^1 - v_j^2 \rangle_{V^*\times V} \leq J_2^0(\gamma (K^k v^1)_{j-1}, \gamma v_j^1; \gamma v_j^2 - \gamma v_j^1 ) +J_2^0(\gamma (K^k v^2)_{j-1}, \gamma v_j^2; \gamma v_j^1 - \gamma v_j^2 ).
\end{align*}
Again $H(A)$(c) and $H(J)$(c) yield
\begin{align*}
  &\big(m_A - m_{J1}c_\gamma^2\big) \|  v_j^1 - v_j^2\|_V^2 \leq 0.
\end{align*}
Under assumption $(H_s)$, we obtain that $v^1_j = v^2_j$. This equality holds for $j=1,\dots,N$, hence a~solution of~$P_{incl}^h$ is unique.

Now, in order to prove $(\ref{e22})$, we set $w^h  =-v_j$ in ($\ref{9}$) to obtain
\[\langle A v_j+B(K^kv)_{j-1},-v_j\rangle_{V^*\times V}+J_2^0(\gamma (K^k v)_{j-1},\gamma v_j;-\gamma v_j
\geq \langle f_j,  - v_j \rangle_{V^*\times V}.\]
Hence,
\begin{equation}
\langle A v_j,  v_j \rangle_{V^*\times V}\leq J_2^0(\gamma (K^k v)_{j-1}, \gamma v_j; -\gamma v_j)
+\langle B (K^kv)_{j-1}, - v_j \rangle_{V^*\times V}  + \langle  f_j , v_j \rangle_{V^*\times V}. \label{c5}
\end{equation}
Using $H(J)$(c), we get
\begin{align*}
J_2^0(\gamma (K^k v)_{j-1}, \gamma v_{j};  - \gamma v_{j} ) + J_2^0(0, 0; \gamma v_{j} )
  &\leq m_{J1}\|  \gamma v_{j}\|_X^2 + m_{J2}\|  \gamma (K^k v)_{j-1}\|_X\|  \gamma v_{j}\|_X\\
  &\leq m_{J1}c_\gamma^2\| v_{j}\|_V^2 + m_{J2}c_\gamma^2\| (K^k v)_{j-1}\|_V\|v_j\|_V.
\end{align*}
From Proposition 3.23(ii) in \cite{MOS} and assumption $H(J)$(b) we have
\[ J_2^0(0, 0; \gamma v_{j} ) \leq c_0\| \gamma v_{j}\|_X \leq c_0 c_\gamma\| v_{j}\|_V.\]
Hence,
\begin{align}
  J_2^0(\gamma (K^k v)_{j-1}, \gamma v_{j};  - \gamma v_{j} ) \leq m_{J1}c_\gamma^2\| v_{j}\|_V^2 + m_{J2}c_\gamma^2\| (K^k v)_{j-1}\|_V\|  v_{j}\|_V + c_0c_\gamma\| v_{j}\|_V.  \label{c4}
\end{align}
From $H(B)$,
\begin{align}
&\langle B (K^kv)_{j-1}, v_j \rangle_{V^*\times V} \leq  \|B(K^k v)_{j-1}\|_{V^*} \| v_{j}\|_V \leq (L_B  \|(K^k v)_{j-1}\|_V +c) \| v_{j}\|_V. \label{c1}
\end{align}
Using $H(A)$(c), (\ref{c4}) and (\ref{c1}) in (\ref{c5}), we get
\[  m_A \|v_j\|_V^2 \leq m_{J1}c_\gamma^2\| v_{j}\|_V^2 + \big(L_B + m_{J2}c_\gamma^2\big)\|(K^k v)_{j-1}\|_V\|v_j\|_V
+ c\|v_j\|_V + \|f_j\|_{V^*}\|v_j\|_V
\]
Then
\[ \big(m_A-m_{J1}c_\gamma^2-\varepsilon\big)\|v_j\|_V^2\le c_\varepsilon\|(K^k v)_{j-1}\|_V^2+(\|f_j\|_{V^*}+c)\,\|v_j\|_V
. \]
Taking sufficiently small $\varepsilon > 0$, using ($H_s$) and inequality \eqref{ele},
we obtain
\[ \|v_j\|_V^2 \leq ck^2N\sum_{i=1}^{j-1} \| v_i\|_V^2 + c\| u_0\|_V^2 + c (\|f_j\|_{V^*} + 1)\,\| v_j\|_V
\]
From discrete version of the Gronwall inequality (Lemma \ref{gronwall}) with
\[ e_j = \|v_j\|_V^2, \qquad
g_j = c\| u_0\|_V^2 + c(\|f_j\|_{V^*} + 1)\,\| v_j\|_V
,\]
we get
\[\|v_j\|_V^2 \leq  c(\|f_j\|_{V^*} + 1)\,\| v_j\|_V, \]
and then
\[ \|v_j\|_V \leq  c(1+\|f_j\|_{V^*}). \]
which concludes the proof of (\ref{e22}).
\end{pfs}

We now consider an optimization problem, which is equivalent to Problem~$P_{incl}^h$ under the stated assumptions.
To this end, let us fix $j \in \{1,\dots,N\}$ and assume $v_0,\dots,v_{j-1}$ to be known, hence also $(K^kv)_{j-1}$
is given. Let a~functional $\mathcal{L}_j\colon V \rightarrow \mathbb{R}$ be defined for all $ w \in V$ as follows
\begin{equation}
\mathcal{L}_j(w)=\frac{1}{2}\langle Aw ,w\rangle_{V^*\times V}
+\langle B (K^kv)_{j-1} -  f_j,  w \rangle_{V^*\times V} + J (\gamma (K^k v)_{j-1}, \gamma w). \label{L}
\end{equation}

The next lemma collects some properties of the functional~$\mathcal{L}_j$.

\begin{lemma} \label{Lprop}
Under assumptions  $H(A)$, $H(B)$, $H(J)$, $H(f)$, $H(u_0)$ and $(H_s)$, the functional $\mathcal{L}_j\colon V \rightarrow \mathbb{R}$ defined by {\rm{(\ref{L})}} has the following properties
\begin{enumerate}[{\rm (i)}]
 \item $\mathcal{L}_j$ is locally Lipschitz continuous,
 \item $\partial \mathcal{L}_j( w ) \subseteq A w + B (K^kv)_{j-1} - f_j + \gamma^* \, \partial_2 J(\gamma (K^k v)_{j-1}, \gamma w)$ for all $w \in V$,
 \item $\mathcal{L}_j$ is strictly convex,
 \item $\mathcal{L}_j$ is coercive.
\end{enumerate}
\end{lemma}
\begin{pfs}
The proof of (i) is immediate since for $v \in V$ the functional $\mathcal{L}_j$ is locally Lipschitz continuous as a sum of locally Lipschitz continuous functions with respect to $w$.

\smallskip
For the proof of (ii), we observe that from $H(A)$, $H(B)$ and $H(f)$, the functions
\begin{align*}
&\psi_1\colon V \ni  w \mapsto \frac{1}{2} \langle A w , w  \rangle_{V^*\times V} \in \mathbb{R}, \quad
\psi_2\colon V \ni w \mapsto \langle  B (K^kv)_{j-1} - f_j, w  \rangle_{V^*\times V} \in \mathbb{R}
\end{align*}
are strictly differentiable with
\begin{align*}
\psi_1'( w ) = Aw,\quad \psi_2'( w ) = B (K^kv)_{j-1} - f_j.
\end{align*}
Now, using the sum and the chain rules for generalized subgradient (cf.\ Propositions~3.35 and 3.37 in $\cite{MOS}$), we obtain
\begin{align*}
\partial \mathcal{L}_j(w ) &= \psi_1'( w )+\psi_2'( w ) + \partial_2 (J\circ \gamma)(\gamma (K^k v)_{j-1}, w)\\
&\subseteq A w + B (K^kv)_{j-1} - f_j + \gamma^* \partial_2 J(\gamma (K^k v)_{j-1}, \gamma w),
\end{align*}
which concludes (ii).

\smallskip
In order to prove (iii), let us fix $ w^1, w^2 \in V$. We take $ \eta_i \in \partial \mathcal{L}_j( w^i)$ for $i=1, 2$.
From (ii), there exists $ \zeta_i \in \partial_2J(\gamma (K^k v)_{j-1} , \gamma w^i)$ such that
\begin{align*}
\eta_i = A w^i + B (K^kv)_{j-1} - f_j + \gamma^* \zeta_i.
\end{align*}
From the equivalent condition to $H(J)$(c), and consequently from~(\ref{relaxed}), we have
\begin{align*}
\langle  \partial_2J(\gamma (K^k v)_{j-1} , \gamma w^1)
 - \partial_2J(\gamma (K^k v)_{j-1} , \gamma w^2), \gamma w^1 - \gamma w^2  \rangle_{X^*\times X}
  \geq -m_{J1}\|\gamma w^1 - \gamma w^2  \|_X^2.
\end{align*}
Hence and from $H(A)$(c), we obtain
\begin{align*}
\langle \eta  _1 -  \eta  _2,  w^1 -  w^2 \rangle_{V^*\times V}
&=\langle A w^1 - A w^2,  w^1-w^2\rangle_{V^*\times V}+\langle \gamma^*\zeta_1-\gamma^*\zeta_2,w^1-w^2\rangle_{V^*\times V}\\
&\geq m_A\| w^1 -  w^2\|_V^2 + \langle \zeta_1 - \zeta_2, \gamma w^1 - \gamma w^2  \rangle_{X^*\times X} \\
&\geq m_A\| w^1 -  w^2\|_V^2  -  m_{J1} \|\gamma w^1 - \gamma w^2  \|_X^2 \\
&\geq (m_A - m_{J1} c_\gamma^2)\| w^1 -  w^2\|_V^2.
\end{align*}
From $(H_s)$ we see that $\partial \mathcal{L}_j$ is strongly monotone. This is equivalent to the fact that $\mathcal{L}_j$ is strongly convex (see Theorem~3.4 in $\cite{FLG}$), which implies that it is strictly convex.

\smallskip
For the proof of (iv), let us fix $w \in V$. From $H(A)$(c) we obtain
\begin{align}
\mathcal{L}_j(w) &=\frac{1}{2} \langle  Aw , w  \rangle_{V^*\times V} +\langle  B (K^kv)_{j-1}, w  \rangle_{V^*\times V} -  \langle  f_j,  w \rangle_{V^*\times V}+ J (\gamma (K^k v)_{j-1}, \gamma w)\nonumber\\
&\ge\frac{1}{2} m_A \|w\|_V^2 - \|B (K^kv)_{j-1}\|_{V^*}\|w\|_V- \|f_j\|_{V^*}\|w\|_V +J (\gamma (K^k v)_{j-1},\gamma w).
  \label{lab1}
\end{align}
Now, using the Lebourg mean value theorem (cf. Proposition~3.36 in \cite{MOS}), we get that there exist $\lambda \in (0,1)$ and $\eta \in \partial_2(J\circ\gamma)(\gamma (K^k v)_{j-1}, \lambda w)$ such that
\begin{equation}
J(\gamma  (K^k v)_{j-1}, \gamma w) = \langle \eta, w  \rangle_{V^*\times V} + J(\gamma  (K^k v)_{j-1}, 0). \label{c2}
\end{equation}
Since $\partial_2(J \circ \gamma)(\gamma (K^k v)_{j-1}, \lambda w) \subseteq \gamma^* \partial_2 J(\gamma (K^k v)_{j-1} , \lambda \gamma w)$ we have $\eta \in \gamma^* \partial_2 J(\gamma (K^k v)_{j-1}, \lambda \gamma w )$. Then there exists $z_1 \in \partial_2 J(\gamma (K^k v)_{j-1} , \lambda \gamma w )$ such that $\eta =  \gamma^* z_1$. Taking $z_2 \in \partial_2 J(\gamma (K^k v)_{j-1}, 0)$ and by (\ref{relaxed}), we obtain
\begin{align*}
  &\lambda\langle \gamma^* z_1 -  \gamma^* z_2, w  \rangle_{V^*\times V} = \langle z_1 - z_2, \lambda \gamma w \rangle_{X^*\times X} \geq - m_{J1} \| \lambda \gamma w \|_X^2 \geq - m_{J1} \lambda^2 c_\gamma^2 \| w \|_V^2,
\end{align*}
and this, along with the fact that $\lambda \in (0,1)$, leads to
\begin{equation}
\langle\eta, w\rangle_{V^*\times V}\ge - m_{J1}c_\gamma^2\|w\|_V^2+\langle z_2,\gamma w\rangle_{X^*\times X}. \label{c30}
\end{equation}
Using $H(J)$(b), we get for given $v$
\begin{align*}
  \|\partial_2 J(\gamma (K^k v)_{j-1}, 0)\|_{X^*} \leq c_0 + c_1 \|\gamma (K^k v)_{j-1}\|_{X} \leq c.
\end{align*}
Hence
\begin{align}
\begin{split}
  \langle z_2, \gamma w \rangle_{X^*\times X} \geq -|\langle z_2, \gamma w \rangle_{X^*\times X}| \geq - \|\partial_2 J(\gamma (K^k v)_{j-1}, 0)\|_{X^*}\| \gamma w \|_X \geq -c\|w \|_V. \label{c40}
\end{split}
\end{align}
Combining (\ref{lab1})--(\ref{c40}) and because $J(\gamma (K^k v)_{j-1}, 0)$ is bounded from below for fixed first argument,  we get
\begin{align*}
\mathcal{L}_j(w) &\geq \frac{1}{2} m_A \|w\|_V^2 - \|B (K^kv)_{j-1}\|_{V^*}\|w\|_V - \|f_j\|_{V^*}\|w\|_V
- m_{J1} c_\gamma^2 \| w \|_V^2 -c\|w \|_V  - c\\
&\geq \Big(\frac{1}{2}m_A - m_{J1} c_\gamma^2\Big) \| w \|_V^2 -c\|w \|_V - c.
\end{align*}
From $(H_s)$ we see that $\mathcal{L}$ is coercive.
\end{pfs}

The problem under consideration reads as follows.

\medskip
\noindent
\textbf{Problem $\bm{P_{opt}^h}$:} Find $v^{hk} = \{v_j^{hk}\}_{j=1}^N \subset V^h$ such that for $j\in\{1,\dots,N\}$,
\begin{align*}
  0 \in & \partial \mathcal{L}_j( v_j^{hk} ).
\end{align*}

We are now in a position to prove the existence and uniqueness result for the above optimization problem.

\begin{lemma} \label{O}
Assume $H(A)$, $H(B)$, $H(J)$, $H(f)$, $H(u_0)$ and $(H_s)$. Then Problem~$P_{opt}^h$ has a unique solution.
\end{lemma}
\begin{pfs}
From Lemma \ref{Lprop} (i), (iv), we see that the functional $\mathcal{L}_j$ is proper, lower semicontinuous
and coercive. This implies that it attains a global minimum. Uniqueness of that
minimum is guaranteed by Lemma \ref{Lprop} (iii).
\end{pfs}

Let us conclude results from the previous lemma by the following theorem.

\begin{theorem} \label{I=O}
Assume $H(A)$, $H(B)$, $H(J)$, $H(f)$, $H(u_0)$ and $(H_s)$. Then Problems~$P_{incl}^h$ and $P_{opt}^h$ are equivalent,
have a unique solution $v^h\in V^h$ and this solution satisfies
\begin{equation}
\| v_j^{hk} \|_V \leq c(1+\| f_j \|_{V^*}),\quad 1\le j\le N \label{e2}
\end{equation}
with a positive constant $c$.
\end{theorem}
\begin{pfs}
Lemma~$\ref{Lprop}$\,(ii) implies that every solution to Problem $P_{opt}^h$ solves Problem~$P_{incl}^h$.
Using this fact, Lemmas $\ref{I}$ and $\ref{O}$, we see that a unique solution to Problem~$P_{opt}^h$ is also a unique solution to Problem~$P_{incl}^h$. Because of the uniqueness of the solution to Problem~$P_{incl}^h$ we get that Problems~$P_{incl}^h$ and $P_{opt}^h$ are equivalent. The estimate in the statement of the theorem follows from Lemma~$\ref{I}$.
\end{pfs}

Now let us conclude this section by presenting the following main theorem concerning error estimate of the numerical scheme.

\begin{theorem} \label{estimate}
Assume that $H(A)$, $H(B)$, $H(J)$, $H(f)$ and $(H_s)$ hold. Let $v$ and $v^{hk}$ be unique solutions to
Problems~$P_{hvi}$ and $P_{opt}^h$, respectively. Additionally assume that $v\in W^{1,\infty}(0,T;V)$.
Then there exists a constant $c>0$ such that
\begin{align}
\max_{1\leq j \leq N} \|v_j-v ^{hk}_j\|_V^2
& \le c\max_{1\le j\le N}\,\inf\limits_{w^h\in V^h}\Big\{k^2+\|v_j-w^h\|_V^2+\|\gamma v_j-\gamma w^h\|_X
+|R_j(v_j-w^h)| \Big\} \nonumber\\
&\quad{} + c\,\|u_0 - u_0^h\|_{V}^2,   \label{thminequality}
\end{align}
where a residual quantity is given by
\begin{equation*} 
R_j(w)=\langle A v_j,w\rangle_{V^*\times V}+\langle B(K v)_j, w\rangle_{V^*\times V}-\langle f_j, w\rangle_{V^*\times V}.
\end{equation*}
\end{theorem}
\begin{pfs}
Let $v$ be a solution to Problem~$P_{hvi}$ and $v^{hk}$ be a solution to Problem~$P_{opt}^h$,
and let us fix  $j \in \{1,\dots,N\}$.  Then we have
\begin{equation}
\langle  Av_j + B (Kv)_j,  w\rangle_{V^*\times V} + J_2^0(\gamma (Kv)_j, \gamma v_j;  \gamma w)
\geq \langle  f_j,  w\rangle_{V^*\times V} \quad \forall\, w \in V,\label{estim1}
\end{equation}
and, as in (\ref{9}),
\begin{equation}
\langle Av^{hk}_j+B(K^k v^{hk})_{j-1},w^h\rangle_{V^*\times V}+J_2^0(\gamma(K^k v^{hk})_{j-1},\gamma v^{hk}_j; \gamma w^h)\geq \langle  f_j,  w^h\rangle_{V^*\times V} \quad \forall\,w^h \in V^h.\label{estim2}
\end{equation}
Taking (\ref{estim1}) with $w=v^{hk}_j-v_j$ and adding to (\ref{estim2}) with $w^h$ replaced by $w^h - v^{hk}_j $,
after some calculation we obtain for all $w^h \in V^h$,
\begin{align}
\begin{split}
  &\langle A v^{hk}_j - A v_j,  v^{hk}_j - v_j \rangle_{V^*\times V} \leq\langle  A v^{hk}_j - A v_j,  w^h - v_j  \rangle_{V^*\times V}\\
  &\quad+\langle B (Kv)_j - B (K^k v^{hk})_{j-1}, v^{hk}_j - v_j \rangle_{V^*\times V}\\
  &\quad+ \langle B (K^k v^{hk})_{j-1} - B (K v)_j, w^h - v_j \rangle_{V^*\times V}  \\
  &\quad+ J_2^0(\gamma (K v)_j, \gamma v_j;  \gamma v^{hk}_j - \gamma v_j) + J_2^0(\gamma (K^k v^{hk})_{j-1}, \gamma v^{hk}_j;  \gamma w^h - \gamma v^{hk}_j)\\
  &\quad+ \langle  A v_j,  w^h - v_j  \rangle_{V^*\times V} + \langle B (K v)_j, w^h - v_j \rangle_{V^*\times V} + \langle f_j, v_j - w^h \rangle_{V^*\times V} \label{estim3}
\end{split}
\end{align}
Using $H(A)$ we have
\begin{align}
  &m_A \| v^{hk}_j - v_j\|_V^2 \leq \langle A v^{hk}_j - A v_j,  v^{hk}_j - v_j \rangle_{V^*\times V}. \label{estim4}
\end{align}
From the Schwartz inequality, the fact that $A\in\mathcal{L}(V,V^*)$ and the Young inequality, we get
\begin{align}
\langle  A v^{hk}_j - A v_j,  w^h - v_j  \rangle_{V^*\times V}
  &\leq \| A v^{hk}_j - A v_j\|_{V^*}\|w^h - v_j\|_V \leq \|A\|_{\mathcal{L}(V,V^*)}\| v^{hk}_j - v_j\|_V\| w^h - v_j\|_V \nonumber \\
  &\leq \varepsilon \| v^{hk}_j - v_j\|_V^2 + c_\varepsilon\, \|A\|^2_{\mathcal{L}(V,V^*)} \| w^h - v_j\|_V^2. \label{estim5}
\end{align}
Analogously for the operator $B$,
\begin{equation}
\langle B (Kv)_j - B (K^k v^{hk})_{j-1}, v^{hk}_j - v_j \rangle_{V^*\times V}
\leq L_B^2 c_\varepsilon \| (K^k v^{hk})_{j-1} - (Kv)_j\|_V^2 + \varepsilon\|v^{hk}_j - v_j\|_V^2, \label{estim6}
\end{equation}
and
\begin{equation}
\langle B (K^k v^{hk})_{j-1} - B (K v)_j, w^h - v_j \rangle_{V^*\times V}
\leq\frac{L_B^2}{2} \| (K^k v^{hk})_{j-1} - (Kv)_j\|_V^2 + \frac{1}{2}\|w^h - v_j\|_V^2. \label{estim7}
\end{equation}
Now, from (\ref{e2}) there exits $M > 0$ such that for all $i \in \{1,\dots,N\}$
\[ \| v_i^{hk} \|_V \leq c(1+\| f_i \|_{V^*}) \leq c(1+\max_{1\leq i \leq N}\| f_i \|_{V^*}) \leq M, \]
and, since $kN = T$, we have
\begin{align*}
   \|(K^k v^{hk})_{j-1}\|_V \leq  k \sum_{i=1}^{j-1} \|v_i^{hk}\|_V  + \|u_0^h\|_V \leq kN M + \|u_0^h\|_V \leq T M + c.
\end{align*}
From these estimations combined with Proposition 3.23(ii) in \cite{MOS} and by assumption $H(J)$(b) we have
\begin{align*}
|J_2^0(\gamma (K^k v^{hk})_{j-1}, \gamma v^{hk}_j;  \gamma w^h - \gamma v_j )|
&\leq(c_0 +  c_1c_\gamma \| (K^k v^{hk})_{j-1}\|_V + c_2c_\gamma\| v^{hk}_j\|_V) \|  \gamma w^h -  \gamma v_j \|_X\\
&\leq c \| \gamma w^h -  \gamma v_j \|_X.
\end{align*}
 We now use assumption $H(J)$(c) and the Young inequality to obtain
 \begin{align}
 \begin{split}
  &J_2^0(\gamma (K^k v^{hk})_{j-1}, \gamma v^{hk}_j;  \gamma v_j - \gamma v^{hk}_j )+ J_2^0(\gamma (K v)_j, \gamma v_j;  \gamma v^{hk}_j - \gamma v_j )\\
  &\quad\leq m_{J1}c_\gamma^2\, \|  v^{hk}_j -  v_j \|_V^2 + m_{J2}^2  c_\varepsilon c_\gamma^4\, \|  (K^k v^{hk})_{j-1} - (K v)_j \|_{V}^2 + \varepsilon \, \|  v^{hk}_j - v_j \|_{V}^2.
  \end{split}
\end{align}
Hence, from subadditivity of generalized directional derivative (Proposition 3.23(i) in \cite{MOS}) we get
\begin{align}
\begin{split}
  &J_2^0(\gamma(K^k v^{hk})_{j-1},\gamma v^{hk}_j;  \gamma w^h - \gamma v^{hk}_j ) + J_2^0(\gamma (K v)_j, \gamma v_j;  \gamma v^{hk}_j - \gamma v_j )\\
  &\quad\leq J_2^0(\gamma (K^k v^{hk})_{j-1}, \gamma v^{hk}_j;  \gamma w^h - \gamma v_j ) + J_2^0(\gamma (K^k v^{hk})_{j-1}, \gamma v^{hk}_j;  \gamma v_j - \gamma v^{hk}_j ) \\
  &\qquad+ J_2^0(\gamma (K v)_j, \gamma v_j;  \gamma v^{hk}_j - \gamma v_j )\\
  &\quad\leq c\, \|  \gamma w^h -  \gamma v_j \|_X + m_{J1}c_\gamma^2\, \|  v^{hk}_j -  v_j \|_V^2\\
  &\qquad+ m_{J2}^2  c_\varepsilon c_\gamma^4\, \|  (K^k v^{hk})_{j-1} - (K v)_j \|_{V}^2 + \varepsilon \, \|  v^{hk}_j - v_j \|_{V}^2. \label{estim9}
  \end{split}
\end{align}
Using inequalities $(\ref{estim4})$--$(\ref{estim9})$ in $(\ref{estim3})$, we obtain
\begin{align}
  &(m_A - m_{J1}c_\gamma^2 - 3\varepsilon) \| v^{hk}_j - v_j\|_V^2\nonumber \\
  &\qquad\leq  c_\varepsilon\| (K^k v^{hk})_{j-1} -  (K v)_j \|_{V}^2 + c_\varepsilon\|v_j - w^h \|_{V}^2 + c\, \|  \gamma w^h -  \gamma v_j \|_X\nonumber\\
  &\qquad\quad{} + \langle  A v_j,  w^h - v_j  \rangle_{V^*\times V} + \langle B (K v)_j, w^h - v_j \rangle_{V^*\times V} + \langle  f_j,  v_j - w^h  \rangle_{V^*\times V}. \label{estim10}
\end{align}
We now calculate
 \begin{align}
 \| (K^k v^{hk})_{j-1} - (K v)_j \|_{V}^2 \leq  2(\| (K^k v)_{j-1} - (K v)_j \|_{V}^2 + \| (K^k v^{hk})_{j-1} - (K^k v)_{j-1} \|_{V}^2 ) \label{estim11},
\end{align}
where
\[ (K^k v)_j = k \sum_{i=1}^j v_i + u_0^h,  \quad v \in C([0,T];V). \]
Using the triangle inequality, the elementary inequality \eqref{ele},
and the Jensen inequality (Lemma~\ref{jensen}) with $I = [t_{i-1},t_i]$ (cf.\ (\ref{JenIneq})), since $kN=T$, we obtain
 \begin{align}
 \begin{split}
&\| (K v)_j - (K^k v)_{j-1}\|_{V}^2 =  \| \int_0^{t_j} v(s) \, ds + u_0 - (k \sum_{i=1}^{j-1} v_i + u_0^h) \|_{V}^2\\ &=  \| \sum_{i=1}^{j-1} \int_{t_{i-1}}^{t_{i}}  (v(s) - v_i) \, ds +  \int_{t_{j-1}}^{t_{j}}  v(s)\, ds  + u_0 - u_0^h\|_{V}^2\\
&\leq \Big(\sum_{i=1}^{j-1} \int_{t_{i-1}}^{t_{i}} \| v(s) - v_i \|_{V} \, ds + \int_{t_{j-1}}^{t_{j}}  \|v(s)\|_V\, ds + \|u_0 - u_0^h\|_{V}\Big)^2\\
&\leq cN \sum_{i=1}^{j-1} \big(\int_{t_{i-1}}^{t_i} \| v(s) - v_i \|_{V} \, ds\big)^2 + c \big(\int_{t_{j-1}}^{t_{j}}  \|v(s)\|_V\, ds\big)^2 + c\|u_0 - u_0^h\|_{V}^2\\
&\leq cT \sum_{i=1}^{j-1} \int_{t_{i-1}}^{t_i} \| v(s) - v_i \|_{V}^2 \, ds + c\big(\int_{t_{j-1}}^{t_{j}}  \|v(s)\|_V\, ds\big)^2 + c\|u_0 - u_0^h\|_{V}^2. \label{estim12}
\end{split}
\end{align}
Since $v \in W^{1,\infty}(0,T;V)$, we obtain for all $i \in \{1,\dots,N\}$, $s \in (t_{i-1},t_i] $
\[ \| v(s) -  v_i \|_{V} = \|\int_{t_i}^{s} v'(\tau)\, d\tau \|_{V} \leq \int_{t_i}^{t_{i+1}} \|v'(\tau)\|_{V} d\tau
\leq k\|v'\|_{L^{\infty}(0,T;V)}\leq k\|v\|_{W^{1,\infty}(0,T;V)},\]
and
\[ \int_{t_{j-1}}^{t_{j}}  \|v(s)\|_V ds \leq k\|v\|_{L^{\infty}(0,T;V)}\leq k\|v\|_{W^{1,\infty}(0,T;V)}. \]
Hence, from the above inequalities, we get
\begin{align}
\| (K v)_j - (K^k v)_{j-1}\|_{V}^2
&\leq ck^2 \sum_{i=1}^{j-1} \int_{t_{i-1}}^{t_i}\|v\|_{W^{1,\infty}(0,T;V)}^2\, ds+ ck^2\|v\|_{W^{1,\infty}(0,T;V)}^2
+ c\|u_0 - u_0^h\|_{V}^2 \nonumber\\
&\leq cTk^2 + ck^2 + c\|u_0 - u_0^h\|_{V}^2. \label{estim15}
\end{align}
And now, from the triangle inequality
\begin{align}
\| (K^k v^{hk})_{j-1} - (K^k v)_{j-1} \|_{V}^2
&= \|k\sum_{i=1}^{j-1} v_i^{hk} + u_0^h - (k\sum_{i=1}^{j-1} v_i + u_0^h)\|_{V}^2\nonumber\\
&\leq \Big(k\sum_{i=1}^{j-1} \|v_i^{hk}- v_i\|_V\Big)^2\leq Tk\sum_{i=1}^{j-1}\| v_i^{hk} - v_i\|_{V}^2. \label{estim16}
\end{align}
Combining (\ref{estim11}), (\ref{estim15}) and (\ref{estim16}), we get
\[ \| (K^k v^{hk})_{j-1} - (K v)_j \|_{V}^2 \leq ck^2 + ck\sum_{i=1}^{j-1}\| v_i^{hk}-v_i\|_V^2 +c\|u_0-u_0^h\|_{V}^2. \]
Returning to (\ref{estim10}), we obtain
\begin{align*}
(m_A - m_{J1}c_\gamma^2 - 3\varepsilon) \| v^{hk}_j - v_j\|_V^2
&\leq c_\varepsilon k^2 + c_\varepsilon\|v_j - w^h \|_{V}^2 + c\, \|\gamma w^h -\gamma v_j\|_X+|R_j(v_j-w^h)|\nonumber \\
&\quad{} + c_\varepsilon\|u_0 - u_0^h\|_{V}^2 + c_\varepsilon k\sum_{i=1}^{j-1}\| v_i^{hk} - v_i\|_{V}^2.
\end{align*}
Taking sufficiently small $\varepsilon$, from smallness assumption $(H_s)$ we get
\begin{align*}
\| v^{hk}_j - v_j\|_V^2 &\leq c k^2 + c\|v_j - w^h \|_{V}^2 + c\, \|  \gamma w^h -  \gamma v_j \|_X \nonumber\\
  &\quad{}+ c |R_j(v_j-w^h)| + c\|u_0 - u_0^h\|_{V}^2 + c k\sum_{i=1}^{j-1}\| v_i^{hk} - v_i\|_{V}^2.
\end{align*}
From discrete version of the Gronwall inequality (Lemma \ref{gronwall}) with
\begin{align*}
& e_j = \| v^{hk}_j - v_j\|_V^2, \\
& g_j = c k^2 + c\|v_j - w^h \|_{V}^2 + c\, \|\gamma w^h -\gamma v_j \|_X + c |R_j(v_j-w^h)| + c\|u_0 - u_0^h\|_{V}^2,
\end{align*}
we obtain the required estimation.
\end{pfs}

\section{Application to a contact problem} \label{sectionACP}

In this section we apply the results of previous sections to a sample mechanical contact problem.
Let us start by introducing the physical setting and notation.

A viscoelastic body occupies domain $\Omega \subset \mathbb{R}^{d}$, where $d = 2,3$ in application.
We assume that its boundary $\Gamma$ is divided into three disjoint measurable parts:
$\Gamma_{D}$, $\Gamma_{C}$, $\Gamma_{N}$, where $\Gamma_{D}$ has a positive measure.
Additionally $\Gamma$ is Lipschitz continuous, and therefore the outward normal vector $\bm{\nu}$ to $\Gamma$ exists
a.e.\ on the boundary. The body is clamped on $\Gamma_{D}$, i.e. its displacement is equal to $\bm{0}$ on this part
of boundary.  A surface force of density $\bm{f}_N$ acts on the boundary~$\Gamma_{N}$ and a body force of density
$\bm{f}_0$ acts in $\Omega$. The contact phenomenon on $\Gamma_{C}$ is modeled using general subdifferential inclusion. Forces and contact conditions may be time dependent. We are interested in finding  body displacement in the time interval $[0,T]$, with $T>0$.

Let us denote by ``$\cdot$'' and $\|\cdot\|$ the scalar product and the Euclidean norm in $\mathbb{R}^{d}$ or $\mathbb{S}^{d}$, respectively, where  $\mathbb{S}^{d} = \mathbb{R}^{d \times d}_{sym}$. Indices $i$ and $j$ run from $1$ to $d$ and the index after a comma represents the partial derivative with respect to the corresponding component of the independent variable.
Summation over repeated indices is implied. We denote the divergence operator by $\textrm{Div }\bm{\sigma} = (\sigma_{ij,j})$. The standard Lebesgue and Sobolev spaces $L^2(\Omega)^d = L^2(\Omega;\mathbb{R}^d)$ and $H^1(\Omega)^d = H^1(\Omega;\mathbb{R}^d)$ are used.
The linearized (small) strain tensor for displacement $\bm{u} \in H^1(\Omega)^d$ is defined by
\begin{equation}\nonumber
 \bm{\varepsilon}(\bm{u})=(\varepsilon_{ij}(\bm{u})), \quad \varepsilon_{ij}(\bm{u}) = \frac{1}{2}(u_{i,j} + u_{j,i}).
\end{equation}
Let $u_\nu= \bm{u}\cdot \bm{\nu}$ and $\sigma_\nu= \bm{\sigma}\bm{\nu} \cdot \bm{\nu}$ be the normal components of $\bm{u}$ and $\bm{\sigma}$, respectively,  and let $\bm{u}_\tau =\bm{u}-u_\nu\bm{\nu}$ and $\bm{\sigma}_\tau =\bm{\sigma}\bm{\nu}-\sigma_\nu\bm{\nu}$ be their tangential components, respectively. In what follows, for simplicity, we sometimes do not indicate explicitly the dependence of various functions on the spatial variable $\bm{x}$. Now let us introduce the classical formulation of considered quasistatic mechanical contact problem.

\medskip
\noindent
\textbf{Problem $\bm{P^M}$:} \textit{Find a displacement field $\bm{u}\colon \Omega \times [0,T] \rightarrow \mathbb{R}^{d}$ and a stress field $\bm{\sigma}\colon \Omega \times [0,T] \rightarrow \mathbb{S}^{d}$ such that for all $t \in[0,T]$},
\begin{align}
\bm{\sigma}= \mathcal{A}(\bm{\varepsilon}(\bm{u}'))+\mathcal{B}(\bm{\varepsilon}(\bm{u}))\quad&\textrm{ in } \Omega, \label{P1}\\
\textrm{Div }\bm{\sigma} + \bm{f_{0}} = \bm{0} \quad &\textrm{ in } \Omega, \label{P2}\\
\bm{u} = \bm{0}  \qquad &\textrm{ on } \Gamma_{D}, \label{P3}\\
\bm{\sigma\nu} = \bm{f_{N}} \qquad &\textrm{ on } \Gamma_{N}, \label{P4}\\
-\sigma_{\nu} = g_\nu(u_{\nu}) \qquad &\textrm{ on } \Gamma_{C}, \label{P5}\\
-\bm{\sigma_{\tau}} \in g_\tau (u_{\nu})\,\partial j_{\tau}(\bm{u_{\tau}}') \qquad &\textrm{ on } \Gamma_{C}, \label{P6}\\
\bm{u}(0) = \bm{u}_{0} \qquad &\textrm{ in } \Omega. \label{P7}
\end{align}

Here, equation~\eqref{P1} represents an viscoelastic constitutive law, $\mathcal{A}$ is a viscosity operator and
$\mathcal{B}$ is an elasticity operator. Equilibrium equation~(\ref{P2}) reflects the fact that problem is quasistatic.
Equation~(\ref{P3}) represents clamped boundary condition on $\Gamma_{D}$ and~(\ref{P4}) reflects the forces acting
on $\Gamma_{N}$. Relation~(\ref{P5}) describes the response of the foundation in normal direction, whereas the friction
is modeled by inclusion~(\ref{P6}), where $j_\tau$ is a given superpotential, and $g_\tau$ is a~given friction bound.
Finally, equation~(\ref{P7}) represents the initial condition with the initial displacement~$\bm{u}_0$.

We use the following Hilbert spaces
\begin{align*}
&H = L^2(\Omega; \mathbb{R}^d), \quad \mathcal{H} = L^2(\Omega; \mathbb{S}^d),\\
 &H_1 = \{\bm{u} \in H\ |\ \bm{\varepsilon}(\bm{u}) \in \mathcal{H}\},
\quad \mathcal{H}_1 = \{\bm{\sigma} \in \mathcal{H}\ |\ \textrm{Div }\bm{\sigma} \in H\},\\
 &V = \{\bm{v} \in H_1\ |\ \bm{v} = 0 \textrm{ on } \Gamma_{D}\}
\end{align*}
endowed with the inner products
\begin{align*}
&(\bm{u},\bm{v})_H=\int_\Omega u_iv_i\,dx,\quad(\bm{\sigma},\bm{\tau})_\mathcal{H}=\int_\Omega \sigma_{ij}\tau_{ij}\,dx, \\
&(\bm{u},\bm{v})_{H_1} = (\bm{u},\bm{v})_H +(\bm{\varepsilon}(\bm{u}),\bm{\varepsilon}(\bm{v}))_\mathcal{H}, \quad
(\bm{\sigma},\bm{\tau})_{\mathcal{H}_1}=(\bm{\sigma},\bm{\tau})_\mathcal{H}
+(\textrm{Div}\,\bm{\sigma},\textrm{Div}\,\bm{\tau})_H,\\
&(\bm{u},\bm{v})_V = (\bm{\varepsilon}(\bm{u}),\bm{\varepsilon}(\bm{v}))_\mathcal{H}
\end{align*}
and corresponding norms $\|\cdotp\|_X$ with $X$ being $H, \mathcal{H}, H_1, \mathcal{H}_1, V$.
The fact that space $V$ equipped with the norm $\|\cdot\|_V$ is complete follows from Korn's inequality, and its application is allowed because we assumed that $meas(\Gamma_{D}) > 0$.
We consider the trace operator $\gamma \colon V \to L^2(\Gamma_{C})^d=X$.
By the Sobolev trace theorem we know that $\gamma \in  \mathcal{L}(V, X)$ with the norm equal to $c_\gamma$.

\smallskip
Now we present the hypotheses on data of Problem~$P^M$.

\smallskip
\noindent
$\underline{H({\mathcal{A}})}:$ \quad ${\mathcal{A}} \colon \Omega \times {\mathbb S}^d \to {\mathbb S}^d$ satisfies
\begin{enumerate}
  \item[(a)]
     $\mathcal{A}(\bm{x},\bm{\tau}) = (a_{ijkh}(\bm{x})\tau_{kh})$
     for all $\bm{\tau} \in {\mathbb S}^d$, a.e. $\bm{x}\in\Omega,\ a_{ijkh} \in L^{\infty}(\Omega),$
  \item[(b)]
    $\mathcal{A}(\bm{x},\bm{\tau}_1) \cdot \bm{\tau}_2 = \bm{\tau}_1 \cdot  \mathcal{A}(\bm{x},\bm{\tau}_2)$ for all $\bm{\tau}_1, \bm{\tau}_2 \in {\mathbb S}^d$, a.e. $\bm{x}\in\Omega$,
  \item[(c)]
     there exists $m_{\mathcal{A}}>0$ such that $\mathcal{A}(\bm{x},\bm{\tau}) \cdot \bm{\tau} \geq m_{\mathcal{A}} \|\bm{\tau}\|^2$ for all $\bm{\tau} \in {\mathbb S}^d$, a.e. $\bm{x}\in\Omega$.
\end{enumerate}

\noindent
$\underline{H({\mathcal{B}})}:$ \quad ${\mathcal{B}} \colon \Omega \times {\mathbb S}^d \to {\mathbb S}^d$ satisfies
\begin{enumerate}
  \item[(a)]
     ${\mathcal B}(\cdot,\bm{\tau})$ is measurable on $\Omega$ for all
        $\bm{\tau}\in \mathbb{S}^d$, ${\mathcal B}(\cdot,\bm{\tau})\in {\mathcal H}$,
  \item[(b)]
    there exists $L_{\cal B} > 0$ s.t.\ $\|{\cal B}(\bx,\bm{\tau}_1)-{\cal B}(\bx,\bm{\tau}_2)\|
        \le L_{\cal B}\|\bm{\tau}_1- \bm{\tau}_2\|$ for all $\bm{\tau}_1,\bm{\tau}_2\in \mathbb{S}^d$, a.e.\ $\bx\in\Omega$,
  \item[(c)]
     there exists $m_{\mathcal{B}}>0$ such that $(\mathcal{B}(\bm{x},\bm{\tau}_1)-\mathcal{B}(\bm{x},\bm{\tau}_2))\cdot (\bm{\tau}_1-\bm{\tau}_2) \geq m_{\mathcal{B}} \|\bm{\tau}_1-\bm{\tau}_2\|^2$ for all $\bm{\tau}_1, \bm{\tau}_2 \in {\mathbb S}^d$, a.e. $\bm{x}\in\Omega$.
\end{enumerate}

\noindent
$\underline{H(j_{\tau})}:$ \quad $j_{\tau} \colon \Gamma_C \times \mathbb{R}^{d} \to \mathbb{R}$ satisfies
\begin{enumerate}
  \item[(a)]
    $j_{\tau}(\cdot, \bm{\xi})$ is measurable on $\Gamma_C$ for all $\bm{\xi} \in \mathbb{R}^{d}$ and there exists $\bm{e} \in L^2(\Gamma_C)^{d}$ such that $j_{\tau}(\cdot,\bm{e}(\cdot))\in L^1(\Gamma_C)$,
  \item[(b)]
    there exists $c_{\tau}>0$ such that \\[2mm]
    \hspace*{1cm}$|j_{\tau}(\bm{x}, \bm{\xi}_1) - j_{\tau}(\bm{x}, \bm{\xi}_2)| \leq c_{\tau} \|\bm{\xi}_1 - \bm{\xi}_2\|$\quad for all $\bm{\xi}_1, \bm{\xi}_2 \in \mathbb{R}^d$, a.e. $\bm{x} \in \Gamma_C$,
  \item[(c)]
    there exists $\alpha_{\tau} \geq 0$ such that  \\[2mm]
   \hspace*{1cm} $j_{\tau}^0(\bm{x},\bm{\xi}_1;\bm{\xi}_2-\bm{\xi}_1) + j_{\tau}^0(\bm{x},\bm{\xi}_2;\bm{\xi}_1-\bm{\xi}_2)\leq \alpha_\tau\|\bm{\xi}_1-\bm{\xi}_2\|^2$\\[2mm]
     for all $\bm{\xi}_1, \bm{\xi}_2 \in \mathbb{R}^{d}$, a.e. $\bm{x} \in \Gamma_C$.
\end{enumerate}

\noindent
$\underline{H(g)}:$ \quad $g_\iota\colon\Gamma_C\times\mathbb{R}\to\mathbb{R}$, $\iota\in\{\nu, \tau\}$, satisfies
\begin{enumerate}
  \item[(a)]
    $g_\iota(\cdot, r)$ is measurable on $\Gamma_C$ for all $r \in \mathbb{R}$,
  \item[(b)]
    there exists $\overline{g_\iota} > 0$ such that $0 \leq g_\iota(\bm{x}, r) \leq \overline{g_\iota}$ for all $r \in \mathbb{R}$, a.e. $\bm{x} \in \Gamma_C$,
  \item[(c)]
    there exists $L_{g_\iota}>0$ such that \\[2mm]
    \hspace*{1cm} $|g_\iota(\bm{x}, r_1) - g_\iota(\bm{x}, r_2)| \leq L_{g_\iota} |r_1 - r_2|$\quad for all $r_1, r_2 \in \mathbb{R}$, a.e. $\bm{x} \in \Gamma_C$.
\end{enumerate}

\noindent
$(\underline{H_0}): \quad \bm{f}_0 \in C([0,T]; L^2(\Omega)^d), \quad \bm{f}_N \in C([0,T]; L^2(\Gamma_N)^d)$.

\smallskip
We remark that condition $H(j_{\tau})$(b) is equivalent to the fact that $j_{\tau}(\bm{x},\cdot)$ is 
Lipschitz continuous and there exists $c_{\tau} \geq 0$ such that $\|\partial j_{\tau}(\bm{x},\bm{\xi})\|\leq c_{\tau}$ for all $\bm{\xi} \in \mathbb{R}^{d}$ and a.e.\ $\bm{x} \in \Gamma_C$. Moreover, condition $H(g)$(b) is enough to obtain presented mathematical results, but from mechanical point of view we should additionally assume that $g_\iota(r)=0$ for $r\le 0$. This corresponds to the situation when body is separated from the foundation and normal response of the foundation and friction force vanishes.

\smallskip
\noindent
Using the standard procedure, the Green formula and the definition of generalized subdifferential, we obtain a weak formulation of Problem~$P^M$ in the form of hemivariational inequality.

\smallskip
\noindent
\textbf{Problem $\bm{P^M_{hvi}}$:} {\it \ Find a velocity $\bm{v}\colon [0,T]\rightarrow V$ satisfying
for all\ $\bm{w}\in V$ and\  $t\in [0,T]$}
\begin{align}
&\langle A\bm{v}(t)+B(K\bm{v})(t),\bm{w}\rangle_{V^*\times V}+\int_{\Gamma_C}j_2^0
(\gamma(K\bm{v})(t),\gamma\bm{v}(t);\gamma\bm{w})\,da\ge\langle\bm{f}(t),\bm{w}\rangle_{V^*\times V}. \label{HVI1}
\end{align}

\smallskip
\noindent
Here, operators $A,B \colon V \to V^*$, $K\colon L^2(0,T;V) \to L^2(0,T;V)$ and $\bm{f}(t) \in V^*$ are defined for all
$\bm{w},\bm{v} \in V$, $t \in [0,T]$ as follows:
\begin{align*}
&\langle A\bm{v},\bm{w}\rangle_{V^*\times V}=(\mathcal{A}(\bm{\varepsilon}(\bm{v})),\bm{\varepsilon}(\bm{w}))_\mathcal{H},\\
&\langle B\bm{v},\bm{w}\rangle_{V^*\times V}=(\mathcal{B}(\bm{\varepsilon}(\bm{v})),\bm{\varepsilon}(\bm{w}))_\mathcal{H},\\
&(K\bm{v})(t) = \int_0^t \bm{v}(t) \, ds + \bm{u}_0,\\
&\langle \bm{f}(t),\bm{w}\rangle_{V^*\times V}=\int_{\Omega}\bm{f}_0(t)\cdot\bm{w}\,dx
+\int_{\Gamma_C}\bm{f_N}(t)\cdot \gamma\bm{w}\, da,
\end{align*}
and $j\colon \Gamma_C\times \mathbb{R}^d \times \mathbb{R}^d \to \mathbb{R}$ be defined for all $\bm{\xi}, \bm{\eta} \in \mathbb{R}^d$ and $\bm{x}\in \Gamma_C$ by
\begin{align}
j(\bm{x}, \bm{\eta}, \bm{\xi}) =  g_{\nu}(\bm{x}, \eta_{\nu})\,\xi_\nu + g_{\tau}(\bm{x}, \eta_{\nu})\, j_{\tau}(\bm{x}, \bm{\xi}_{\tau}). \label{Jj}
\end{align}
It is easy to check that under assumptions $H(\mathcal{A})$, $H(\mathcal{B})$ and $(H_0)$, the operators~$A$ and $B$ and the functional~$\bm{f}$ satisfy $H(A)$, $H(B)$ and $H(f)$, respectively.
We also define the functional $J \colon L^2(\Gamma_C)^d \times L^2(\Gamma_C)^d \to \mathbb{R}$ for all $\bm{w}, \bm{v} \in L^2(\Gamma_C)^d $ by
\begin{align}
J(\bm{w}, \bm{v}) =  \int_{\Gamma_C} j(\bm{x}, \bm{w}(\bm{x}), \bm{v}(\bm{x}))\, da.\label{J}
\end{align}

\medskip
\noindent
Below we present some properties of the functional~$J$.

\begin{lemma} \label{LJ}
Assumptions $H(j_{\tau})$ and $H(h)$ imply that functional $J$ defined by {\rm{(\ref{Jj})--(\ref{J})}} satisfies $H(J)$.
\end{lemma}

\begin{pfs}
We first observe that from $H(j_\tau)$(a),(b) and $H(h)$(a),(c) the function $j(\cdot,\bm{\eta},\bm{\xi})$ is measurable on $\Gamma_C$, there exist $\bm{e}_1, \bm{e}_2 \in L^2(\Gamma_C)^d$ such that $j(\cdot, \bm{e}_1(\cdot), \bm{e}_2(\cdot))\in L^1(\Gamma_C)$, $j(\bm{x}, \cdot, \bm{\xi})$ is continuous and $j(\bm{x}, \bm{\eta}, \cdot)$ is locally Lipschitz. Moreover, by $H(j_\tau)$(b) and $H(h)$(b) we easily conclude
\begin{align*}
\| \partial_2 j(\bm{x}, \bm{\eta},\bm{\xi})\| &\leq g_\nu(\bm{x},\eta_\nu) + g_\tau(\bm{x},\eta_\nu) \|\partial j_\tau (\bm{x},\bm{\xi}_\tau)\|
\leq \overline{g_\nu}+\overline{g_\tau} \, c_\tau \|\bm{\xi}\|.
\end{align*}
Applying similar procedure to one presented in the proof of Corollary~4.15 in $\cite{MOS}$, we obtain that functional $J$ is well defined, locally Lipschitz with respect to the second variable and the growth condition $H(J)$(b) holds with $c_{0}= \sqrt{2 meas(\Gamma_C)}\, \overline{g_\nu}$, $c_{1}=\sqrt{2}\, \overline{g_\tau}\, c_\tau$ and  $c_{2}=0$.

\medskip
\noindent
To prove $H(J)$(c), we take $\bm{\eta}_i, \bm{\xi}_i \in \mathbb{R}^d$, $i = 1, 2$, and by the sum rules
(cf.\ Proposition~3.35 in \cite{MOS}) and from $H(j_\tau)$(b),(c) and $H(h)$(b),(c), we obtain
\begin{align*}
& j_2^0(\bm{x},\bm{\eta}_1,\bm{\xi}_1;\bm{\xi}_2-\bm{\xi}_1)+j_2^0(\bm{x},\bm{\eta}_2,\bm{\xi}_2;\bm{\xi}_1- \bm{\xi}_2)\\
&\qquad \leq g_{\nu}(\bm{x},\eta_{1\nu})(\xi_{2\nu}- \xi_{1\nu}) - g_{\nu}(\bm{x},\eta_{2\nu})(\xi_{1\nu}- \xi_{2\nu})  \\
&\qquad\quad{} + g_{\tau}(\bm{x},\eta_{1\nu})\left(j_{\tau}^0(\bm{x},\bm{\xi}_{1\tau};
\bm{\xi}_{2\tau}-\bm{\xi}_{1\tau}) + j_{\tau}^0(\bm{x},\bm{\xi}_{2\tau}; \bm{\xi}_{1\tau}-\bm{\xi}_{2\tau})\right) \\
&\qquad\quad{} + \big(g_{\tau}(\bm{x}, \eta_{2\nu})-g_{\tau}( \bm{x},\eta_{1\nu})\big)\,j_{\tau}^0(\bm{x},\bm{\xi}_{2\tau}; \bm{\xi}_{1\tau} -  \bm{\xi}_{2\tau}) \\
&\qquad \leq \overline{g_\tau}\,\alpha_\tau \, \|\bm{\xi}_1 - \bm{\xi}_2\|^2
  + (L_{g_\nu} + L_{g_\tau} c_\tau) \|\bm{\eta}_1 - \bm{\eta}_2\|\, \|\bm{\xi}_1 - \bm{\xi}_2\|.
\end{align*}
Consequently, since
\begin{align*}
J_2^0(\bm{w}, \bm{v}; \bm{z}) \leq \int_{\Gamma_C} j_2^0 (\bm{x}, \bm{w}(\bm{x}), \bm{v}(\bm{x}); \bm{z}(\bm{x}))\,da
\end{align*}
(cf.\ Corollary~4.15~(iii) in $\cite{MOS}$), we have
\begin{align*}
&J_2^0(\bm{w}_1, \bm{v}_1; \bm{v}_2- \bm{v}_1) + J_2^0(\bm{w}_2, \bm{v}_2; \bm{v}_1- \bm{v}_2) \\
&\qquad \leq \int_{\Gamma_{C}}\hspace{-4mm} \left(\overline{g_{\tau}} \alpha_{\tau}\|\bm{v}_1(\bm{x})-\bm{v}_2(\bm{x})\|^2
+(L_{g_\nu}+L_{g_\tau}c_{\tau})\|\bm{w}_1(\bm{x})-\bm{w}_2(\bm{x})\|\,\|\bm{v}_1(\bm{x})-\bm{v}_2(\bm{x})\|\right)\, da.
\end{align*}
Hence, by the H\"{o}lder inequality, we obtain $H(J)$(c) with $m_\alpha = \overline{g_{\tau}} \alpha_{\tau}$ and $m_L =L_{g_\nu} + L_{g_\tau}c_{\tau}$.
\end{pfs}

\smallskip
With the above properties, we can check that assumptions of Theorems \ref{existence}, \ref{I=O} and \ref{estimate}
are satisfied. We can employ previously presented abstract framework and conclude that Problem $P^M_{hvi}$
has a~unique solution.

\smallskip
We now turn to the numerical solution of $P^M_{hvi}$.  For simplicity, we assume $\Omega$ is a polygonal/polyhedral
domain, and express the three parts of the boundary, $\Gamma_D$, $\Gamma_N$ and $\Gamma_C$  as unions of closed
flat components with disjoint interiors:
\[ \Gamma_D =\cup_{i=1}^{i_D}\Gamma_{D,i},\quad \Gamma_N =\cup_{i=1}^{i_N}\Gamma_{N,i},
\quad \Gamma_C =\cup_{i=1}^{i_C}\Gamma_{C,i}. \]
Let $\{{\cal T}^h\}_h$ be a regular family of finite element partitions of $\overline{\Omega}$ into
triangular/tetrahedral elements, compatible with the partition of the boundary $\partial\Omega$ into
$\Gamma_{D,i}$ for $1\le i\le i_D$, $\Gamma_{N,i}$ for $1\le i\le i_N$, and $\Gamma_{C,i}$ for $1\le i\le i_C$,
i.e.\ if the intersection of one side/face of an element with one of these sets has a positive measure,
then the side/face lies entirely in that set.  Here $h\to0$ denotes the finite element mesh-size.
Corresponding to the partition ${\cal T}^h$, we introduce the linear finite element space
\begin{equation*}
V^h = \left\{\bv^h\in C(\overline{\Omega})^d \mid \bv^h|_T\in \mathbb{P}_1(T)^d
   \ \forall\,T\in {\cal T}^h,\,\bv^h=\bzero\ {\rm on\ }\Gamma_D\right\}.
\label{NA1}
\end{equation*}
As in the previous section, given a positive integer $N$, define the step-size $k=T/N$ and the nodes $t_n=nk$,
$0\le n\le N$.   We will assume
\begin{equation}
\bm{u}_0\in H^2(\Omega)^d.
\label{4.12b}
\end{equation}
Let $\bm{u}^h_{0}\in V^h$ be the interpolant of $\bm{u}_{0}$ in $V^h$.  Then (\cite{BS2008, Ci1978})
\begin{equation}
\|\bm{u}_0-\bm{u}^h_0\|_V\le c\,h.
\label{4.12a}
\end{equation}

Then we introduce the following discretized version of Problem $P_{hvi}$.

\medskip
\noindent
\textbf{Problem $\bm{P_{hvi}^h}$:} {\it \ Find a velocity $\bm{v}^{hk}=\{\bm{v}^{hk}_j\}_{j=1}^N\subset V^h$
such that for $1\le j\le N$,}
\begin{align*}
&\langle A\bm{v}^{hk}_j+B(K^k \bm{v}^{hk})_{j-1},\bm{w}^h\rangle_{V^*\times V}+\int_{\Gamma_C}j_2^0
(\gamma(K^k \bm{v}^{hk})_{j-1},\gamma\bm{v}^{hk}_j;\gamma\bm{w}^h)\,da\ge\langle\bm{f}_j,\bm{w}^h\rangle_{V^*\times V}
\quad\forall\,\bm{w}^h\in V^h.
\end{align*}

Similar to Problem $P^M_{hvi}$, under the stated assumptions, Problem $P_{hvi}^h$ has a unique solution.
Assume $\bv\in W^{1,\infty}(0,T;V)$.  It is easy to see from the proof of Theorem \ref{estimate} that the
inequality \eqref{thminequality} remains valid when $J^0_2(\bm{u},\bm{v};\bm{w})$ is replaced by
$\int_{\Gamma_C} j^0_2(\bm{u},\bm{v};\bm{w})\,da$.  Thus, we have a constant $c>0$ such that
\begin{align}
\max_{1\leq j \leq N} \|\bv_j-\bv ^{hk}_j\|_V^2
& \le c\max_{1\le j\le N}\,\inf\limits_{{\boldsymbol w}^h\in V^h}\Big\{k^2+\|\bv_j-\bw^h\|_V^2
+\|\gamma \bv_j-\gamma \bw^h\|_{L^2(\Gamma_C)^d}+|R_j(\bv_j-\bw^h)| \Big\} \nonumber\\
&\quad{} + c \|\bu_0 - \bu_0^h\|_{V}^2,   \label{NA11}
\end{align}
where
\begin{equation} \label{NA12}
R_j(\bw)=\langle A\bv_j+ B(K\bv)_j,\bw\rangle_{V^*\times V}-\langle\fb_j,\bw\rangle_{V^*\times V}.
\end{equation}

\begin{theorem}
Assume $H(\mathcal{A})$, $H(\mathcal{B})$, $H(j_{\tau})$, $H(g)$, $(H_0)$, $(H_s)$, and \eqref{4.12b},
and assume the solution regularity
\begin{align}
& \bv\in W^{1,\infty}(0,T;V),\quad \bv\in C([0,T]; H^2(\Omega)^d),\quad \bsigma\bnu\in C([0,T]; L^2(\Gamma_C)^d),
\label{4.14a}\\
& \bv|_{\Gamma_{C,i}}\in C([0,T]; H^2(\Gamma_{C,i})^d),\ 1\le i\le i_C. \label{4.14b}
\end{align}
Then, for the solution $\bm{v}$ to Problem~$P^M_{hvi}$ and the solution $\bm{v}^{hk}$ to Problem~$P_{hvi}^h$ there exists a constant $c>0$ such that
\begin{equation} \label{NA13}
\max_{1\leq j \leq N} \|\bv_j-\bv ^{hk}_j\|_V\leq c\left(k+h\right).
\end{equation}
\end{theorem}
\begin{pfs}
We bound the residual term defined by (\ref{NA12}) using similar procedure to one described in \cite{HS}.
Let $\bm{w}\in C^\infty(\overline\Omega)^d$ be arbitrary with $\bm{w}=\bm{0}$ on $\Gamma_D\cup\Gamma_C$ in
the inequality \eqref{HVI1}. Let us fix $j \in \{1,\dots,N\}$. We can derive the equality
\[\langle A\bm{v}_j+B(K\bm{v})_j, \bm{w} \rangle_{V^*\times V}  = \langle \bm{f}_j, \bm{w} \rangle_{V^*\times V}
\quad \forall\,\bm{w}\in C^\infty(\overline\Omega)^d, \ \bm{w}=\bm{0}\ {\rm on}\ \Gamma_D\cup\Gamma_C.\]
Then it is possible to deduce that
\[ \textrm{Div}\,(\mathcal{A}(\bm{\varepsilon}(\bm{u}'_j)) + \mathcal{B}(\bm{\varepsilon}(\bm{u}_j)))
    + (\bm{f}_{0})_j = \bm{0} \quad{\rm in\ the\ sense\ of\ distributions.} \]
Since $(\bm{f}_{0})_j \in L^2(\Omega)^d$, we have
\begin{equation}
    \textrm{Div}\,(\mathcal{A}(\bm{\varepsilon}(\bm{u}'_j)) + \mathcal{B}(\bm{\varepsilon}(\bm{u}_j)))
    + (\bm{f}_0)_j = \bm{0}  \quad\textrm{a.e.\ in } \Omega. \label{d1}
\end{equation}
It is also possible to deduce that
 \begin{equation}
    \bm{\sigma}_j\bm{\nu} = (\bm{f}_{N})_j \quad \textrm{a.e.\ on } \Gamma_{N}. \label{d2}
\end{equation}
For details, see \cite{HS}.  We multiply equation (\ref{d1}) by $\bm{w}\in V$ to obtain
\[  \int_\Gamma\bm{\sigma}_j\bm{\nu}\cdot\bm{w}\,da-\int_{\Omega}\left[\mathcal{A}(\bm{\varepsilon}(\bm{u}'_j))
+ \mathcal{B}(\bm{\varepsilon}(\bm{u}_j))\right]\cdot\bm{\varepsilon}(\bm{w})\,dx
+\int_\Omega(\bm{f}_0)_j\cdot\bm{w}\,dx=0.\]
Using the homogeneous Dirichlet boundary condition of $\bm{w}$ on $\Gamma_D$ and the traction boundary condition
given by (\ref{d2}) we have
\[ \int_{\Omega}\left[\mathcal{A}(\bm{\varepsilon}(\bm{u}'_j))
+ \mathcal{B}(\bm{\varepsilon}(\bm{u}_j))\right]\cdot\bm{\varepsilon}(\bm{w})\,dx
=\int_{\Gamma_C}\bsigma_j\bnu\cdot\bw\,da + \int_\Omega(\bm{f}_0)_j\cdot\bm{w}\,dx+\int_{\Gamma_N}(\bm{f}_N)_j\cdot\bm{w}\,da.\]
Thus, for $R_j(\bw)$ defined by \eqref{NA12}, we have
\begin{equation}
 R_j(\bw)=  \int_{\Gamma_C} \bm{\sigma}_j\bm{\nu}\cdot\bm{w}\,da\le c\,\|\bm{w}\|_{L^2(\Gamma_C)^d}. \label{R2}
\end{equation}
We denote by $\Pi^h \bm{u} \in V^h$ the finite element interpolant of $\bm{u}$.
From \eqref{NA11}, (\ref{R2}) and \eqref{4.12a},  we get
\begin{equation}
\max_{1\leq j \leq N} \|\bv_j-\bv ^{hk}_j\|_V^2
\le c\max_{1\le j\le N}\left\{k^2+\|\bv_j-\Pi^h \bv_j\|_V^2
+\|\bv_j-\Pi^h \bv_j\|_{L^2(\Gamma_C)^d}\right\}+ c\,h^2.
\label{4.24a}
\end{equation}
By the standard finite element interpolation error bounds (\cite{BS2008, Ci1978}), due to the solution regularity
\eqref{4.14a} and \eqref{4.14b}, we have
\begin{align*}
  &\|\bv_j-\Pi^h \bv_j\|_V \leq c\,h,\\
  &\|\bv_j-\Pi^h \bv_j\|_{L^2(\Gamma_C)^d} \leq c\,h^2.
\end{align*}
Using these bounds in \eqref{4.24a}, we obtain the error estimate \eqref{NA13}.
\end{pfs}

\section{Numerical results} \label{sectionNR}

In this section, we report computer simulation results on a numerical example. We apply numerical scheme~$P_{opt}^h$ to approximate solution of $P^M_{hvi}$. The linear finite element space $V^h$ based on uniform triangular partition of $\overline{\Omega}$ and the uniform partition of the time interval $[0,1]$ with the time step size $k=1/N$ for a positive integer $N$ are used. In order to minimize not necessarily differentiable functional $\mathcal{L}$ we use Powell's conjugate direction method. This method does not require the assumption that optimized function is differentiable. Other, more refined nonsmooth optimization algorithms described for instance in \cite{BKM}, could also be adapted.

\begin{figure}[ht]
\centering
    \includegraphics[width=0.4\linewidth]{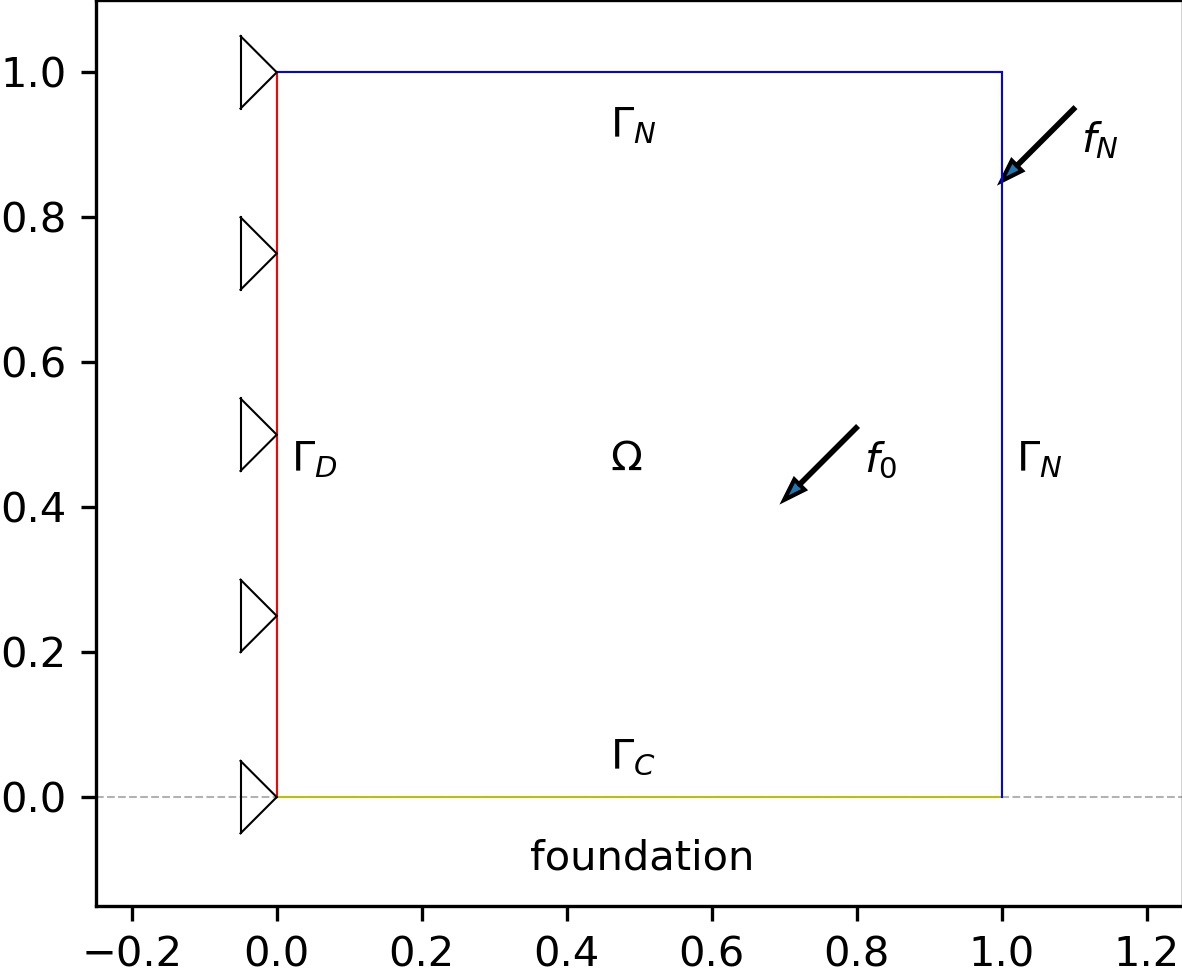}
 \caption{Initial setting} \label{setting}
\end{figure}

We set $d=2$ and consider a square-shaped set $\Omega=(0,1)\times(0,1)$ presented in Figure \ref{setting} with the following partition of the boundary
\[ \Gamma_{D} = \{0\} \times [0,1], \quad \Gamma_{N} = ([0,1] \times \{1\}) \cup (\{1\} \times [0,1]),
\quad \Gamma_{C} = [0,1] \times \{0\}. \]

\noindent We employ the Kelvin-Voigt type short memory viscoelastic law for the isotropic body. The viscosity operator $\mathcal{A}$ and the elasticity operator $\mathcal{B}$ are defined by
\begin{align*}
\mathcal{A}(\bm{\tau}) = 2\phi\bm{\tau} + \xi \mbox{tr}(\bm{\tau})I,\qquad \bm{\tau} \in \mathbb{S}^2,\\
\mathcal{B}(\bm{\tau}) = 2\eta\bm{\tau} + \lambda \mbox{tr}(\bm{\tau})I,\qquad \bm{\tau} \in \mathbb{S}^2.
\end{align*}

\noindent Here $I$ denotes the identity matrix, $\mbox{tr}$ denotes the trace of the matrix, $\lambda$ and $\eta$
are the Lam\'e  coefficients, whereas $\phi$ and $\xi$ represent the viscosity coefficients,  $\lambda, \eta, \phi, \xi > 0$.
In our simulations, we choose $\phi = \xi = 2$, $\lambda = \eta = 4$, $T=1$.

We first demonstrate the effect of some input data on the deformation of the body.
In all cases, we show the shape of the body at final time $t=1$, as well as the contact interface forces on $\Gamma_C$. The numerical solutions correspond
to the time step size $1/32$ and the boundary $\Gamma_C$ of the body divided into $32$ equal parts.
Let us take the following data
\begin{align*}
 &\bm{u}_{0}(\bm{x}) = (0,0), \quad \bm{x} \in \Omega,\\
 &g_\nu(\bm{x}, \eta)= \left \{ \begin{array}{ll}
   0, &\eta \in (-\infty,\, 0), \\
   30\, \eta, &\eta \in [0,\,0.1), \\
   3, &\eta \in [0.1,\infty), \\
  \end{array} \right. \bm{x} \in \Gamma_C,\\
 &g_\tau(\bm{x}, \eta)= g_\nu(\bm{x}, \eta), \quad \eta \in \mathbb{R}, \, \bm{x} \in \Gamma_C,\\
 &j_{\tau}(\bm{x}, \bm{\xi}) = -0.3\,e^{-\|\bm{\xi}\|} + 0.7\,\|\bm{\xi}\| , \quad \bm{\xi} \in \mathbb{R}^2,\ \bm{x} \in \Gamma_C,\\
 &\bm{f}_N(\bm{x},t) = (0,0), \quad \bm{x} \in \Omega,\ t \in [0,T],\\
 &\bm{f}_0(\bm{x},t) = (-2.5,-0.5), \quad \bm{x} \in \Omega,\ t \in [0,T].
\end{align*}

\noindent We note that function $j_{\tau}$, based on Example 7.26 in \cite{MOS}, is nondifferentiable and nonconvex.
Our aim is to investigate reaction of the body to various modifications of input data.

\begin{figure}[ht]
\centering
\begin{minipage}{.45\textwidth}
  \centering
    \includegraphics[width=0.9\linewidth]{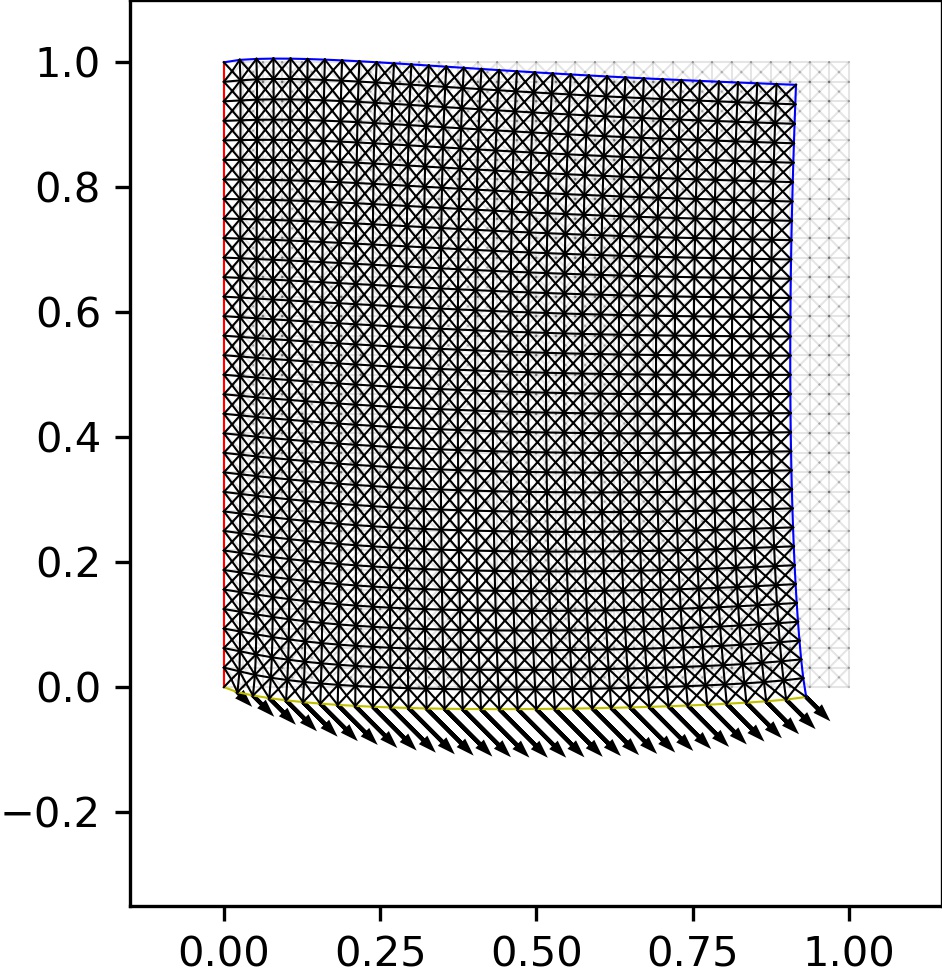}
    \caption{Initial data} \label{figOne}
\end{minipage}
\begin{minipage}{.45\textwidth}
  \centering
    \includegraphics[width=0.9\linewidth]{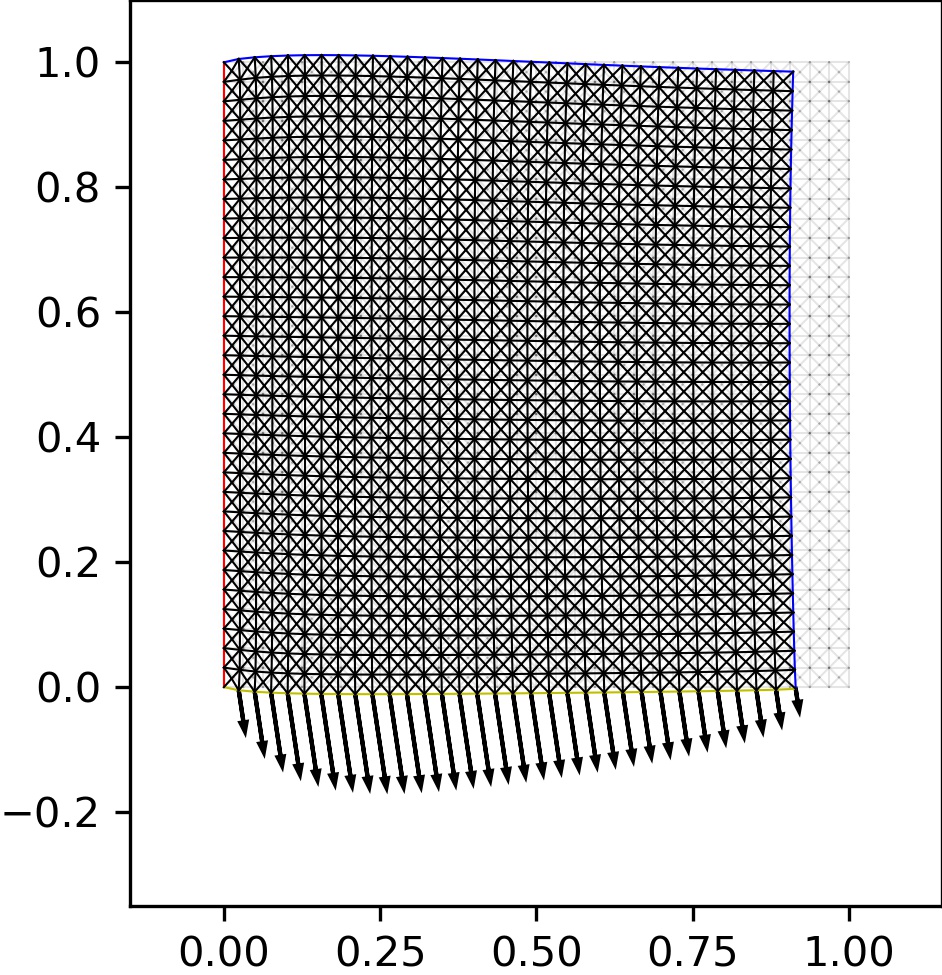}
    \caption{Modified function $g_\nu$} \label{figTwo}
\end{minipage}
\end{figure}

\begin{figure}[ht]
\centering
\begin{minipage}{.45\textwidth}
  \centering
 \includegraphics[width=0.9\linewidth]{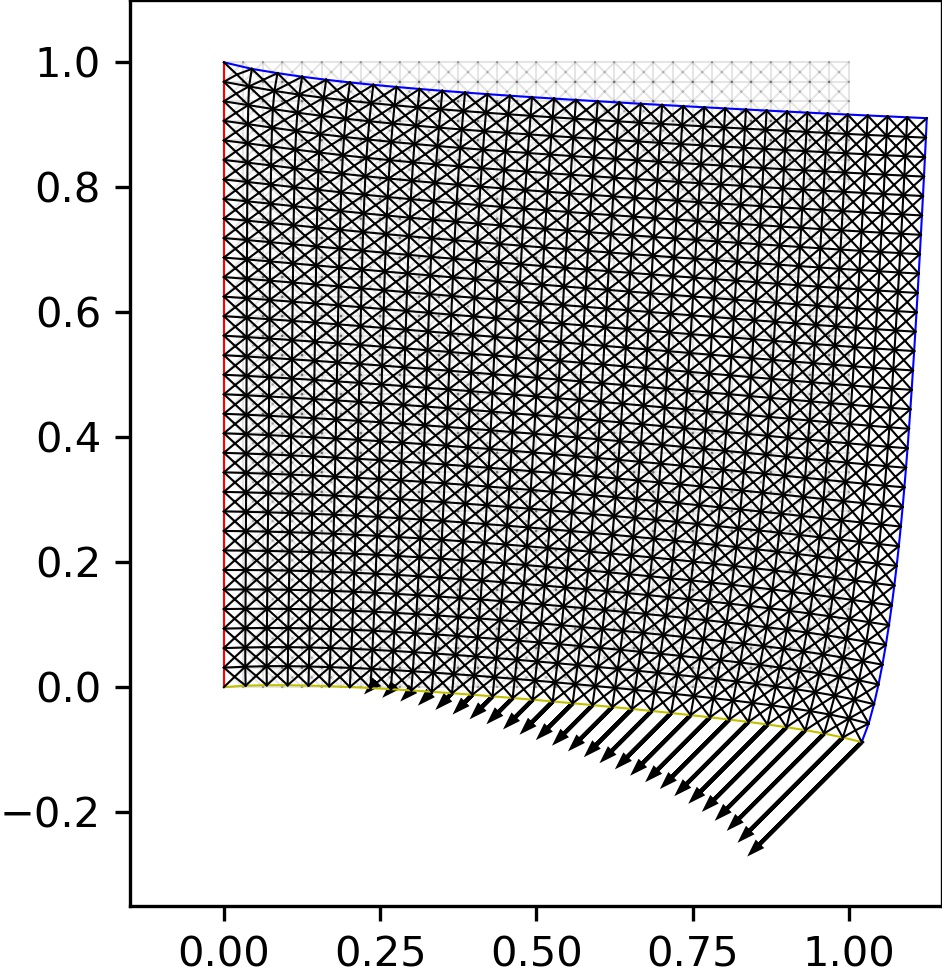}
    \caption{Modified force $\bm{f}_0$} \label{figThree}
\end{minipage}
\begin{minipage}{.45\textwidth}
  \centering
    \includegraphics[width=0.9\linewidth]{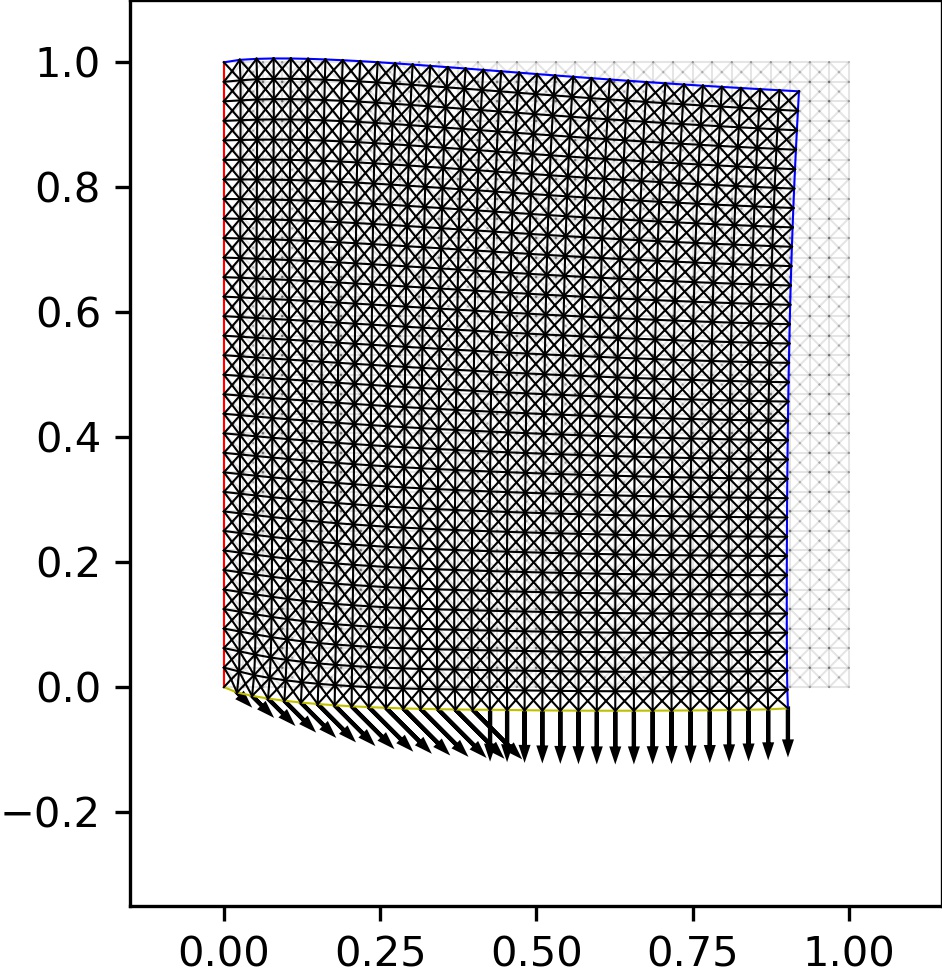}
    \caption{Modified function $g_\tau$} \label{figFour}
\end{minipage}
\end{figure}

In Figure \ref{figOne} we present output obtained without any modifications. We push the body down and to the left
with force $\bm{f}_0$. As a result the body penetrates the foundation, but frictional forces restrict its movement
in proximity of $\Gamma_C$. Next, we modify the function $g_\nu$ to be given by
\begin{align*}
 &g_{\nu}(\bm{x}, \eta) = \left \{ \begin{array}{ll}
   0, &\eta \in (-\infty, 0), \\
   200\, \eta, &\eta \in [0,0.1), \\
   20, &\eta \in [0.1,\infty), \\
  \end{array} \right. \bm{x} \in \Gamma_C.
\end{align*}
In Figure~\ref{figTwo} we observe that this modification models more rigid foundation by increasing its response in normal direction on $\Gamma_C$. The result is decreased penetration of the foundation $u_\nu$. The friction also decreases, due to influence of function $g_\tau$ which depends on $u_\nu$.
Now we return to original data and only change the direction of force $\bm{f}_0$ to the following
\begin{align*}
 &\bm{f}_0(\bm{x},t) = (2.5,\, -0.5), \quad \bm{x} \in \Omega,\ t \in [0,T].
\end{align*}
In Figure~\ref{figThree} we observe that the body displaces in opposite direction to previous examples.
Because of frictional forces it moves to the right more in the higher part than in the lower part.
The penetration of the foundation and friction increase as we get closer to bottom right corner of the body.
In the last experiment we once more return to original data and modify the function $g_\tau$ as follows
\begin{align*}
 &g_{\tau}(\bm{x}, \eta) = \left \{ \begin{array}{ll}
   0, &\eta \in (-\infty, 0),\, \bm{x} \in \Gamma_C,\\
   0, &\eta \in [0,\infty),\, \bm{x} \in [0.5, 1] \times \{0\},\\
   30\, \eta, &\eta \in [0,0.1),\, \bm{x} \in [0, 0.5) \times \{0\}, \\
   3, &\eta \in [0.1,\infty),\, \bm{x} \in [0, 0.5) \times \{0\}.\\
  \end{array} \right.
\end{align*}
In this case $\Gamma_C$ is divided into two parts, and right part is covered in grease. In Figure~\ref{figFour} we see that contact of left part of the body with the foundation creates friction, whereas contact of right part is frictionless.

\begin{table}[ht]
\footnotesize
\centering
\begin{tabular}{ l r r r r r r }\hline
$k$ & $1/2$ & $1/4$ & $1/8$ & $1/16$ & $1/32$ & $1/64$  \\ \hline
$\|\bm{v}-\bm{v}^{hk}\|_V/\|\bm{v}\|_V$ & 3.3088 & 6.6785e--1& 2.1124e--1 & 5.2534e--2&1.4992e--2 & 6.0133e--3\\
{\rm Convergence order} & & 2.3087 & 1.6606 & 2.0076 & 1.8090 & 1.3180\\ \hline
\end{tabular}
\caption{Numerical errors for fixed $h = 1/256$} \label{tabOne}
\end{table}

\begin{table}[ht]
\footnotesize
\centering
\begin{tabular}{ l r r r r r r }\hline
$h$ & $1/2$ & $1/4$ & $1/8$ & $1/16$ & $1/32$ & $1/64$  \\ \hline
$\|\bm{v}-\bm{v}^{hk}\|_V/\|\bm{v}\|_V$ & 2.4390 & 1.4329e--1  & 8.3185e--2 & 4.7945e--2 & 2.7101e--2 & 1.4753e--2\\
{\rm Convergence order} & & 0.7673 & 0.7845 & 0.7949 & 0.8230 & 0.8773\\ \hline
\end{tabular}
\caption{Numerical errors for fixed $k = 1/256$} \label{tabTwo}
\end{table}

\begin{figure}[ht]
\centering
\begin{minipage}{.49\textwidth}
  \centering
  \includegraphics[width=0.9\linewidth]{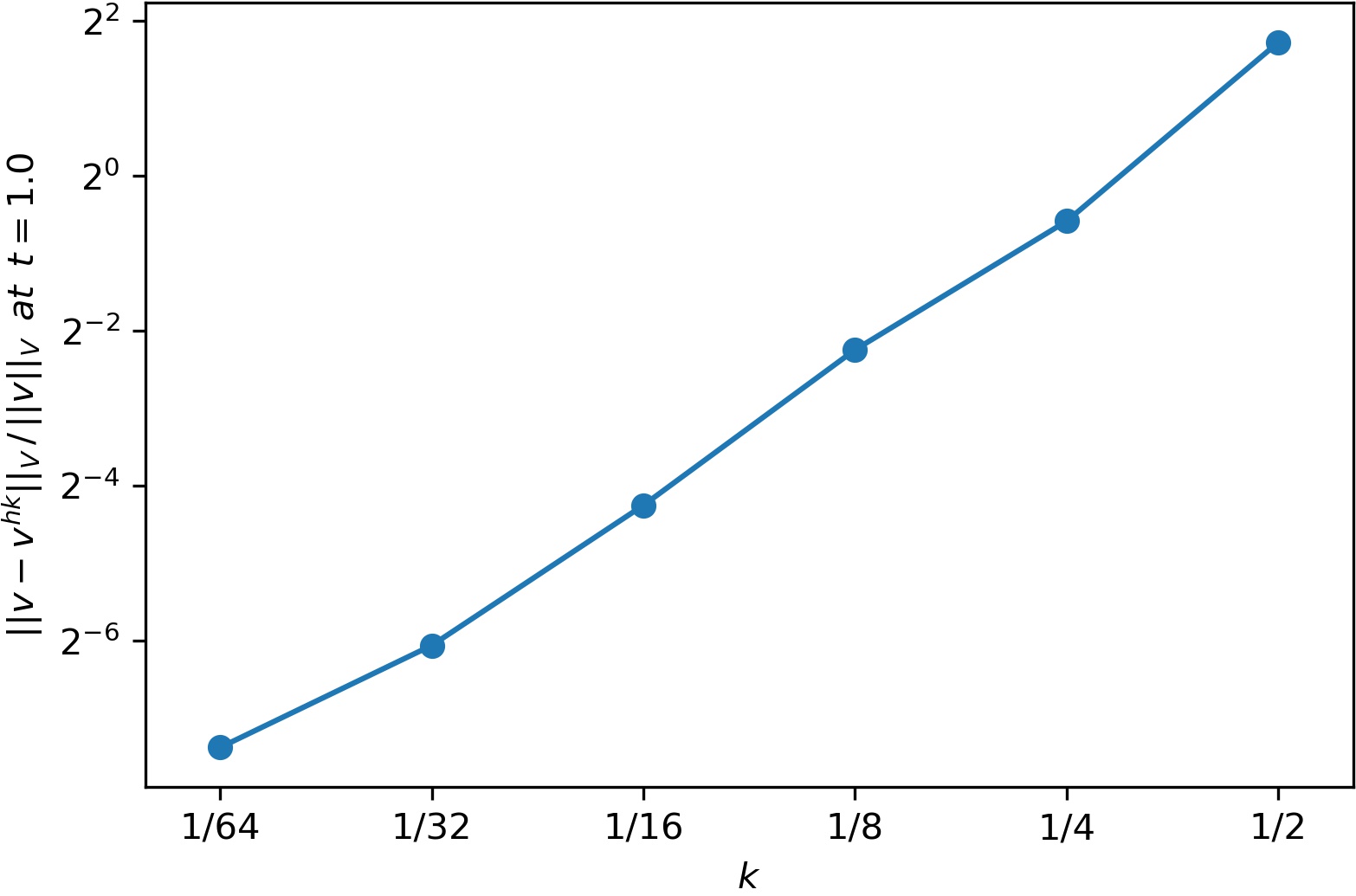}
  \caption{Error estimate, fixed~\mbox{$h = 1/256$}} \label{figFive}
\end{minipage}
\begin{minipage}{.49\textwidth}
  \centering
  \includegraphics[width=0.9\linewidth]{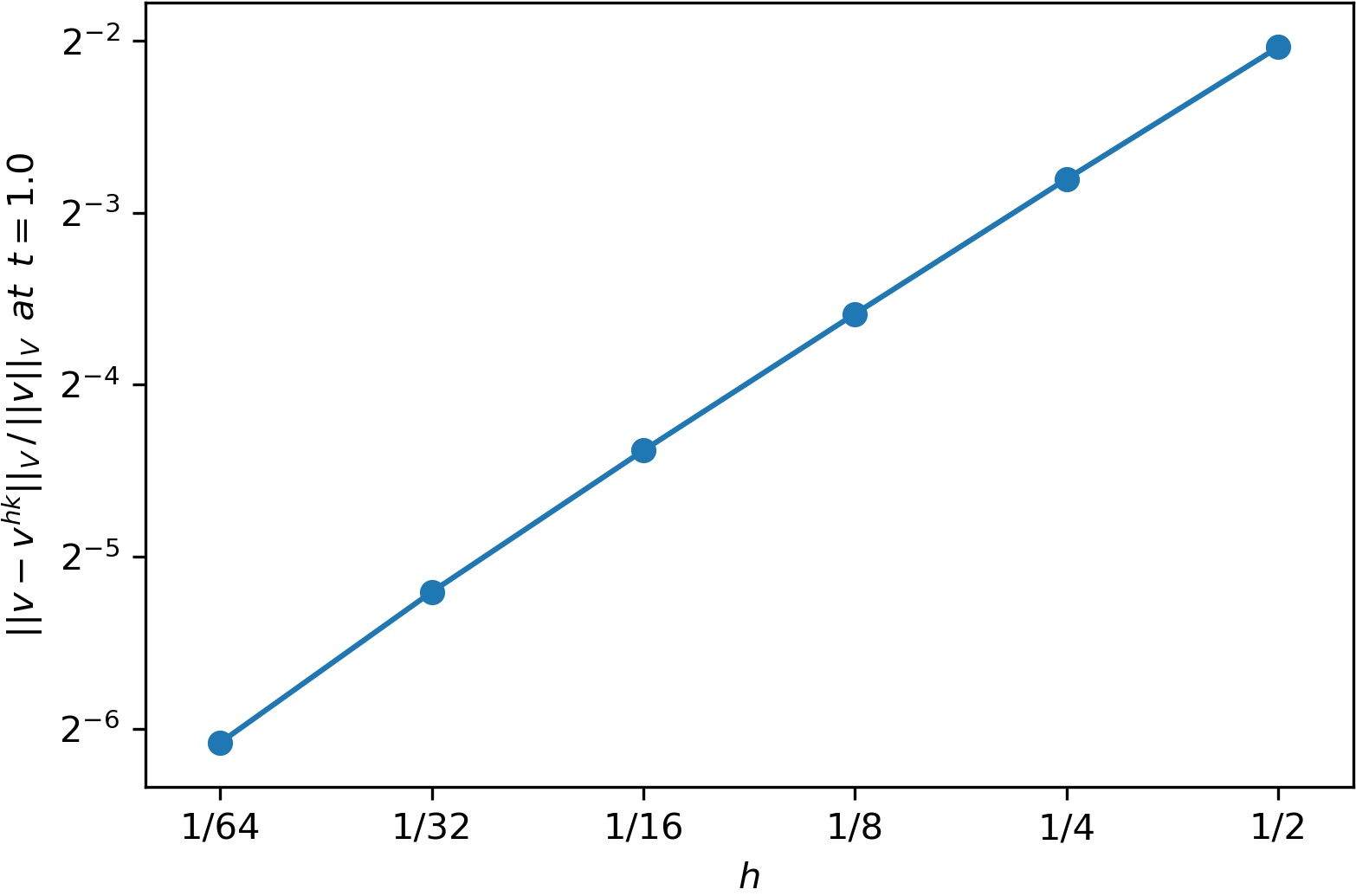}
  \caption{Error estimate, fixed~\mbox{$k = 1/256$}} \label{figSix}
\end{minipage}
\end{figure}

In order to illustrate the error estimate obtained in Section \ref{sectionNS}, we present the empirical convergence orders of the numerical method on the model problem. We take the following data
\begin{align*}
 &\bm{u}_{0}(\bm{x}) = (0,0), \quad \bm{x} \in \Omega,\\
 &g_\nu(\bm{x}, \eta)= \left \{ \begin{array}{ll}
   0, &\eta \in (-\infty,\, 0), \\
   60\, \eta, &\eta \in [0,\,0.1), \\
   6, &\eta \in [0.1,\infty), \\
  \end{array} \right. \bm{x} \in \Gamma_C,\\
 &g_\tau(\bm{x}, \eta)= \left \{ \begin{array}{ll}
   0, &\eta \in (-\infty,\, 0), \\
   120\, \eta, &\eta \in [0,\,0.1), \\
   12, &\eta \in [0.1,\infty), \\
  \end{array} \right. \bm{x} \in \Gamma_C,\\
 &j_{\tau}(\bm{x}, \bm{\xi}) = \|\bm{\xi}\|, \quad \bm{\xi} \in \mathbb{R}^2,\ \bm{x} \in \Gamma_C,\\
 &\bm{f}_N(\bm{x},t) = (-0.2,-0.2), \quad \bm{x} \in \Omega,\ t \in [0,T],\\
 &\bm{f}_0(\bm{x},t) = (-1,-0.4), \quad \bm{x} \in \Omega, \ t \in [0,T].
\end{align*}
We have $\|\bm{v}\|_{V}\doteq 0.06738$.
We present a comparison of numerical errors $\|\bm{v} - \bm{v}^{hk}\|_V$ computed for a sequence of solutions to discretized problems.
We use a uniform discretization of the problem domain and time interval according to the spatial
discretization parameter~$h$ and the time step size~$k$, respectively. The boundary $\Gamma_C$ of $\Omega$
is divided into $1/h$ equal parts. The numerical solution corresponding to $h = 1/256$ and $k = 1/256$ is taken as the ``exact'' solution~$\bm{v}$.

For the first experiment, we fix $h = 1/256$ and start with $k = 1/2$, which is successively halved.
The results are presented in Table \ref{tabOne} and Figure \ref{figFive}, where the dependence of the relative error estimates $\|\bm{v}  - \bm{v}^{hk}\|_V / \|\bm{v}\|_V$  with respect to $k$ are plotted on a log-log scale.  A first order convergence can be observed for the numerical solutions of the displacement.

For the second experiment, we fix $k = 1/256$ and start with $h = 1/2$, which is also successively halved.
The results are presented in Table \ref{tabTwo} and Figure \ref{figSix}. Again, a first order convergence can be observed.

\medskip
\noindent {\bf Acknowledgments}\\
The project has received funding from the European Union's Horizon 2020 Research and Innovation Programme under the Marie Sklo\-do\-wska-Curie grant agreement no.\ 823731 CONMECH.
It is supported by the projects financed by the Ministry of Science and Higher Education of Republic of Poland under
Grants Nos. 4004/GGPJII/H2020/2018/0 and 440328/PnH2/2019.

\end{document}